\documentclass[11pt,letterpaper]{amsart}
\usepackage{amsmath,amsthm,amssymb,upgreek,tipa}
\usepackage{mathabx}
\usepackage[latin1]{inputenc}
\usepackage[T1]{fontenc}
\usepackage[all]{xy}
\usepackage[colorlinks=true,linkcolor=blue,pagebackref,breaklinks]{hyperref}
\usepackage[usenames]{color}
\usepackage{ulem}
\usepackage{url}
\usepackage{marginnote}
\usepackage{enumitem}

\newtheorem{thm}[equation]{Theorem}
\newtheorem*{thm*}{Theorem}

\newtheorem{cor}[equation]{Corollary}
\newtheorem{lem}[equation]{Lemma}
\newtheorem{prop}[equation]{Proposition}
\newtheorem{conj}[equation]{Conjecture}
\newtheorem{hyp}[equation]{Hypothesis}

\theoremstyle{definition}

\newtheorem{rem}[equation]{Remark}
\newtheorem{defn}[equation]{Definition}
\newtheorem{defn-lem}[equation]{Definition-Lemma}
\newtheorem{rmk}[equation]{Remark}
\newtheorem*{rem*}{Remark}
\newtheorem{obsv}[equation]{Observation}
\newtheorem{sub}[equation]{ }

\DeclareMathOperator{\h}{H}

\def\R{\mathbb R}
\def\N{\mathbb N}
\def\Z{\mathbb Z}
\def\A{\mathbb A}
\def\Q{\mathbb Q}
\def\C{\mathbb C}
\def\CC{\mathbb CC}
\def\RR{\mathbb R}
\def\QQ{\mathbb Q}
\def\Qbar{\overline{\Q}}
\def\CC{\mathbb C}
\def\Gm{{\mathbb G}_m}

\def\EE{\mathcal E}
\def\cF{\mathcal F}

\def\ira{\stackrel{\sim}{\longrightarrow}}
\def\hra{\hookrightarrow}
\def\ra{\rightarrow}

\def\g{\mathfrak g}

\def\q{\mathfrak q}

\def\O{\mathcal O}
\def\h{\mathfrak h}

\def\k{\mathfrak k}
\def\p{\mathfrak p}

\def\<{\langle}
\def\>{\rangle}

\def\W{\mathcal W}

\def\Acm{\mathbf{A}}

\newcommand{\Ss}{{\mathbb S}}

\def\Hom{{\rm Hom}}
\def\GL{{\rm GL}}

\def\CL{\mathcal{L}}
\def\CA{\mathcal{A}}

\def\CP{\mathcal{P}}
\def\CO{\mathcal{O}}
\def\cm{F}
\def\tr{F^{+}}
\def\Acm{\mathbb{A}_{F}}

\def\Atr{\mathbb{A}_{F^{+}}}
\def\Atrf{\mathbb{A}_{F^{+},f}}
\newcommand{\isoarrow}{{~\overset\sim\longrightarrow}}

\def\HC{A} 

\begin{document}

\title[]{Factorization of periods, construction of automorphic motives and Deligne's conjecture over CM-fields}
\author{Harald Grobner, Michael Harris \& Jie Lin}
\thanks{H.G. has been supported by the START-prize Y-966 of the Austrian Science Fund (FWF), the FWF Stand-alone research project P32333 and the FWF Principal investigator project PAT-4628923. M.H.'s research received funding from the European Research Council under the European Community's Seventh Framework Program (FP7/2007-2013) / ERC Grant agreement no. 290766 (AAMOT). M.H. was partially supported by NSF Grants DMS-1404769, DMS-1701651, and DMS-2302208.  This work was also supported by the National Science Foundation under Grant No. DMS-1440140 while M.H. was in residence at the Mathematical Sciences Research Institute in Berkeley, California, during the Spring 2019 semester. J.L. was supported by the European Research Council under the European Community's Seventh Framework Programme (FP7/2007-2013) / ERC Grant agreement no. 290766 (AAMOT)}

\maketitle

\begin{abstract}
The present paper is devoted to the relations between Deligne's conjecture on critical values of motivic $L$-functions and the multiplicative relations between periods of arithmetically normalized automorphic forms on unitary groups.  As an application of our main result, we establish Deligne's conjecture for a class of CM-automorphic motives, which we construct in this paper. Our proof uses the results of our recent joint work with Raghuram in combination with the Ichino--Ikeda--Neal-Harris (IINH) formula for unitary groups -- which is now a theorem -- and an analysis of cup products of coherent cohomological automorphic forms on Shimura varieties to establish relations between certain automorphic periods and critical values of Rankin-Selberg and Asai $L$-functions of $\GL(n)\times\GL(m)$ over CM fields.  By reinterpreting these critical values in terms of automorphic periods of holomorphic automorphic forms on unitary groups, we show that the automorphic periods of holomorphic forms can be factored as products of coherent cohomological forms, compatibly with a motivic factorization predicted by the Tate conjecture.  All of these results are stated under a certain regularity condition and an hypothesis of rationality on archimedean zeta-integrals.
\end{abstract}

\setcounter{tocdepth}{1}
\tableofcontents

\section*{Introduction}
Deligne's conjecture, as stated in \cite{deligne}, asserts that the critical values of $L$-functions attached to pure motives $M$ can be expressed in terms of a certain period $c^\pm(M)$ and a power of $(2\pi i)$. This paper is devoted to showing this claim for tensor products of certain CM-motives, which we shall construct in the body of this paper. Let $F$ be a CM field.  For $n \geq 1$ we abbreviate  $G_n:=\GL_n/{\cm}$.  Our first main ingredient to do so is our upcoming paper with Raghuram, in which we are going to prove the following statement:

\begin{thm}[\cite{GHLR,jie-thesis}] \label{automorphic Deligne general intro} 
Let $n,n'\geq 1$ be integers and let $\Pi$ (resp. $\Pi'$) be a cohomological conjugate self-dual cuspidal automorphic representation of $G_n(\Acm)$ (resp. $G_{n'}(\Acm)$), which descends to a tempered cuspidal automorphic representation of a unitary group $U_I(\A_{F^+})$ for each possible $[F^+:\Q]$-tuple of signatures $I$ at the archimedean places, i.e., it satisfies Hypotheis \ref{descent}. We assume that both $\Pi_\infty$ and $\Pi'_{\infty}$ are $5$-regular.   If $n \equiv n' \mod 2$, we assume in addition that the isobaric sum $(\Pi\eta^n)\boxplus (\Pi'^{c}\eta^{n'})$ is $2$-regular; if $n$ and $n'$ have opposite parities then we assume $(\Pi\eta^n)\boxplus (\Pi'^{c}\eta^{n'})$ is $5$-regular. Then the automorphic version of Deligne's conjecture, cf.\ Conjecture \ref{main conjecture}, is true: If $s_0$ is critical, in Deligne's sense, for $L(s,\Pi\times \Pi')$, then the value at $s_0$ of the partial $L$-function $L^S(s,\Pi\times \Pi')$ (for some appropriate finite set $S$), satisfies
\begin{equation}\label{eq:Thm1l}
L^S(s_{0},\Pi\otimes \Pi') \sim_{E(\Pi)E(\Pi')} (2\pi i)^{nn's_{0}} \prod\limits_{\imath \in \Sigma}[\prod\limits_{0\leq i\leq n}P^{(i)}(\Pi,\imath)^{sp(i,\Pi;\Pi',\imath)}\prod\limits_{0\leq j\leq n'}P^{(j)}(\Pi',\imath)^{sp(j,\Pi';\Pi,\imath)}].
\end{equation}
Here $\imath$ runs over complex embeddings of $F$ belonging to a fixed CM type $\Sigma$, $sp(i,\Pi;\Pi',\imath)$ are integers depending on the relative positions of the infinitesimal characters of $\Pi$ and $\Pi'$ at the place $\imath$, and  $P^{(i)}(\Pi,\imath)$ and  $P^{(j)}(\Pi',\imath)$ are period invariants attached to $\Pi$ and $\Pi'$ by quadratic base change from certain unitary groups (that depend on $\imath$ and the superscripts $(i), (j)$), and the symbol ``$\sim_{E(\Pi)E(\Pi')}$'' means that the left-hand side is the product of the right hand side by an element of a certain number field attached to $\Pi$ and $\Pi'$.
\end{thm}
\noindent{\bf Remark.}  The proof in \cite{GHLR}, which is independent of the parity of $n - n'$, requires the isobaric sum $(\Pi\eta^n)\boxplus (\Pi'^{c}\eta^{n'})$ to be $5$-regular.  The proof in \cite{jie-thesis} is valid under the weaker condition given here when $n \equiv n' \mod 2$. \\\\ 
\noindent In other words, the critical values of the Rankin-Selberg $L$-function $L(s,\Pi \times \Pi')$ can be expressed in terms of Petersson norms of certain arithmetically normalized holomorphic automorphic forms.  We refer to \cite{GHLR} and the body of the present paper for details and explanations.  Hypotheis \ref{descent}, a mild local restriction at non-archimedean places, is only relevant when $nn'$ is even, and then amounts to the familiar fact that it is not always possible to construct even-dimensional hermitian spaces with arbitrary local invariants; it can be relaxed by a standard base change construction, as in \cite{yos95}, at the cost of introducing additional quadratic irrationalities.\\\\
Theorem \ref{automorphic Deligne general intro} can be viewed as an {\it automorphic version of Deligne's conjecture} on critical values of motivic $L$-functions, in the special case of the Rankin-Selberg tensor product.  The original conjecture of Deligne states that the left-hand side of the equation in Theorem \ref{automorphic Deligne general intro}  is proportional (up to the coefficient field $E(M(\Pi))E(M(\Pi'))$) to the period invariant Deligne assigned to a motive $R_{F/\Q}(M(\Pi)\otimes M(\Pi'))$ whose $L$-function is given by $L(s,\Pi\otimes \Pi')$; we denote this invariant  $c^+(s_0, R_{F/\Q}(M(\Pi)\otimes M(\Pi')))$.   The main theme of this paper is the relation of the right-hand side of \eqref{eq:Thm1l} to Deligne's period invariant.  Under the hypotheses of Theorem \ref{automorphic Deligne general intro}, cohomological realizations of the hypothetical motives $M(\Pi)$ and $M(\Pi')$ over $F$, of rank $n$ and $n'$ over their respective coefficient fields, can be constructed in the cohomology of Shimura varieties $Sh(V)$ and $Sh(V')$ attached to the unitary groups of hermitian spaces $V$ and $V'$ of rank $n$ and $n'$ respectively.  These Shimura varieties have the property that their connected components are arithmetic quotients of the unit ball in $\C^{n-1}$ and $\C^{n'-1}$, respectively.  The main theorem of the present paper can be paraphrased as follows:

 \begin{thm}\label{cplusintro}   Let $F$, $\Pi$, and $\Pi'$ be as in Theorem \ \ref{automorphic Deligne general intro} and assume in addition that $\Pi$ (resp.\ $\Pi'$) is $(n+4)$-regular (resp.\ $(n+3)$-regular). Then, for any critical value $s_0$ of $L(s,\Pi\times \Pi')$, the Deligne period $c^+(s_0,R_{F/\Q}(M(\Pi)\otimes M(\Pi')))$ can be identified with the right-hand side of the equation in Theorem \ref{automorphic Deligne general intro}, up to a constant that depends only on the infinitesimal characters of $\Pi_\infty$ and $\Pi'_\infty$.   
 \end{thm}
 
\noindent Conjecture\ \ref{lvarch}, which is a statement about the rationality of certain archimedean integrals, asserts that the constant in the last sentence can be taken to be a rational number.  The main theorem is stated in the text as Theorem \ref{main factorization}, which assumes  Conjecture\ \ref{lvarch} as a hypothesis. As with Theorem \ref{automorphic Deligne general intro}, using quadratic base change, as in \cite{yos95}, one can probably obtain a weaker version of Theorem \ref{main factorization} in the absence of this assumption, but we have not checked the details.\\\\ 
The content of Theorem \ref{cplusintro} is a relation between periods of automorphic forms on Shimura varieties attached to hermitian spaces with different signatures. Following an approach pioneered by Shimura over 40 years ago, we combine special cases of Deligne's conjecture with comparisons of distinct expressions for critical values of automorphic $L$-functions to relate automorphic periods on different groups.  
These periods are attached to motives (for absolute Hodge cycles, see \cite[\S 0.9]{deligne}) that occur in the cohomology of the various Shimura varieties.  In view of Tate's conjecture on cycle classes in $\ell$-adic cohomology, the relations obtained are consistent with the determination of the representations of Galois groups of appropriate number fields on the $\ell$-adic  cohomology of the respective motives.  
The paper \cite{linfactorization} used  arguments of this type to show how to factor automorphic periods on Shimura varieties attached to a CM field $F$ as products of automorphic periods of holomorphic modular forms, each attached to an  embedding $\imath = \imath_{v_0}:  F \hookrightarrow \C$.  Theorem \ref{auto-facto} (see \eqref{eq:factor} below) 
leads to a factorization of the latter periods in terms of periods of coherent cohomology classes on Shimura varieties attached to the unitary group $H^{(0)}$ of a hermitian space over $F$ with signature $(n-1,1)$ at $\imath_{v_0}$ and definite at embeddings that are distinct from $\imath_{v_0}$ and its complex conjugate.  
This factorization -- see Theorem \ref{main factorization} (and the explanations in \S\ref{sect:goals}), which is the precise statement of which Theorem \ref{cplusintro} is a paraphrase --  also depends on  Conjecture\ \ref{lvarch} and the local restrictions mentioned above.  The archimedean components $\Pi_\infty, \Pi'_\infty$ of the automorphic representations we consider are tempered and are {\it cohomological}, in the sense that the relative Lie algebra cohomology spaces
\begin{equation}\label{rlag}
H^*(\g_{n,\infty},K_{G_n,\infty},\Pi_\infty\otimes \EE) \neq 0,\quad\quad  H^*(\g_{n-1,\infty},K_{G_{n-1},\infty},\Pi'_\infty\otimes \EE') \neq 0
\end{equation}
for certain finite-dimensional representations $\EE$, $\EE'$ of $G_n$ and $G_{n-1}$, respectively; here $K_{G_n,\infty}$ and $K_{G_{n-1},\infty}$ denote maximal compact mod center subgroups of $G_{n,\infty}$ and $G_{n-1,\infty}$. The local restrictions at archimedean places take the form of regularity hypotheses on the infinitesimal characters of these finite-dimensional representations, or equivalently on the Hodge structures of the associated motives.\\\\
Taken together,  Theorem\ \ref{automorphic Deligne general intro} and Theorem\ \ref{cplusintro} provide a plausible version of Deligne's conjecture for the $L$-function of the tensor products of the motives $M(\Pi)$ and $M(\Pi')$ over $\imath(F)$ attached to $\Pi$ and $\Pi'$, when $\Pi$ and $\Pi'$ satisfy the local restrictions and regularity hypotheses already mentioned -- 
up to the constant that is the subject of Conjecture \ref{lvarch}.  We refer to our Theorem\ \ref{thm:clozel}, where we recall the construction of motives $M(\Pi)$ and $M(\Pi')$. As in \cite[\S 1]{harrisANT}, \cite[\S 5.3]{guer-lin}, what we actually study are collections of realizations in the cohomology of Shimura varieties with coefficients in local systems.  However, since these realizations are obtained as eigenspaces for the action of Hecke correspondences and the local systems are of geometric origin, it is not difficult to see that they can be interpreted as the cohomological realizations attached to Grothendieck motives, using the constructions in \cite[\S III.2]{HT}, for example.  Thus our constructions are compatible with a conjecture of Clozel \cite[Conjecture 4.5]{clozel}. \\\\
 Conjecture \ref{lvarch} can only be settled by a computation of the integrals in question.  The conjecture is natural because its failure would contradict  the Tate conjecture; it is also known to be true in the few cases where it can be checked.  Methods are known for computing these integrals but they are not simple.   In the absence of this conjecture, the methods of this paper provide the weaker statement of Theorem\ \ref{cplusintro}.  A similar statement when $F$ is an imaginary quadratic field had already been proved in \cite{H07} using the theta correspondence, but the proof there is much more complicated.

\subsection*{About the proofs}  
Our main theorem obtains the factorization of periods \eqref{eq:factor} by applying the IINH conjecture to the results on special values, and by using a result on non-vanishing of cup products of coherent cohomology proved in \cite{H14}.  In fact, the case used here had already been treated in \cite{HLi}, assuming properties of stable base change from unitary groups to general linear groups that can be found in \cite{KMSW, {agikms}}. \\\\
The results of \cite{H14} are applied by induction on $n$, and each stage of the induction imposes an additional regularity condition.
This explains the regularity hypothesis in the statement of Theorem \ref{auto-facto}.  The factorization in the theorem must be true in general, but it is not clear to us whether the method based in the IINH conjecture can be adapted in the absence of the regularity hypothesis.

\subsection*{On using the IINH conjecture to solve for unknowns}
Although we have no sympathy with the general outlook of the politician Donald Rumsfeld, and we consider his role in recent history to be largely deleterious, in the formulation of the strategy for proving our main results we did find it helpful to meditate on his thoughts on knowledge, as expressed in the following quotation \cite{EM}:
\\
\\
{\it...as we know, there are known knowns; there are things we know we know. We also know there are known unknowns; that is to say we know there are some things we do not know. But there are also unknown unknowns -- the ones we don't know we don't know.   }
\\
\\
Rumsfeld neglected the {\bf unknown knowns}, such as the period invariants and critical values that are the main subject of this paper.
The formula of Ichino--Ikeda--Neal-Harris, in the inhomogenous form in which it is presented in Theorem \ref{conjecture II}, can be viewed as an identity involving three kinds of transcendental quantities:  critical values of Rankin-Selberg and Asai $L$-functions, Petersson norms of algebraically normalized coherent cohomology classes, and cup products between two such classes.   Here is a simplified version of the conjecture, which is now a theorem, with elementary terms indicated by $(*)$:
\begin{equation}\label{IIsimple}
\cfrac{|I^{can}(f,f')|^{2}}{\<f,f\>\ \<f',f'\>} = (*)\frac{L(\tfrac{1}{2},\Pi\otimes \Pi')}{L(1,\Pi,{\rm As}^{(-1)^n})L(1,\Pi',{\rm As}^{(-1)^{n-1}})}.
\end{equation}
From the Rumsfeld perspective, the denominator of the right-hand side of \eqref{IIsimple}, which is independent of the relative position of $\Pi$ and $\Pi'$,
was an {\bf unknown known} that became a {\bf known known} thanks to \cite{grob-harr} and subsequent generalizations.   The same paper,
as well as \cite{harrisANT}, turn the numerator of the right-hand side into a {\bf known known}, as long as the coefficients $\EE$ and $\EE'$ of the 
cohomology classes defined by $\Pi$ and $\Pi'$ (see \eqref{rlag}) satisfy the relation \eqref{piano}:  
\begin{equation}\label{piano} \Hom_{G_{n-1}}(\EE \otimes \EE', \C) \neq 0. 
\end{equation}

\noindent Theorem \ref{automorphic Deligne general intro}, inspired by a generalization of a formula of Chen \cite{C23}, turns the numerator of the right-hand side into a {\bf known known} even when $\Pi$ and $\Pi'$ themselves do not satisfy \eqref{piano}.\\\\  
Thus the entire right-hand side of the formula can be considered a {\bf known known}.  As for the left-hand side, the periods in the denominator should at best be viewed as {\bf known unknowns}, and then only when $f$ and $f'$ are holomorphic automorphic forms -- because the only thing we know about Petersson norms of (arithmetically normalized) holomorphic automorphic forms is that they are uniquely determined real numbers that are probably transcendental.  That leaves the numerator of the left-hand side, and here we use the result of \cite{H14}, when it applies, to choose $f$ and $f'$ so that the numerator, as a cup product in coherent cohomology, belongs to a fixed algebraic number field.  In fact, the numerator can be taken to be $1$, which is a {\bf  known known}, if anything is.\\\\  
Finally, as the unitary groups vary most of the periods that appear in the numerator of the left-hand side of \eqref{IIsimple} have no cohomological interpretation. Thus these have to be viewed as {\bf unknown unknowns} in Rumsfeld's sense -- precisely because the identity \eqref{IIsimple} relates these periods to {\bf known knowns} and {\bf unknown knowns} (the latter when the periods in the denominator are attached to higher coherent cohomology classes on Shimura varieties for unitary groups with mixed signature, which we have not studied).   We nevertheless believe that the methods of this paper will provide an expression of these {\bf unknown unknowns} in terms of the {\bf known unknowns} which are Petersson norms of arithmetically normalized coherent cohomology classes.   Specifically to $\Pi$ and to a given complex embedding $\imath$ of $F$ we define an invariant $P_i(\Pi,\imath)$, which is the Petersson norm of an automorphic form that defines a(n arithmetically normalized) coherent cohomology class on a Shimura variety attached to a unitary groups definite at all  places in $\Sigma$ other than $\imath$ and of signature $(1,n-1)$ at $\imath$; the index $i$ refers to the degree in coherent cohomology.  The main result of this paper is to show that the Petersson norms of (arithmetically normalized) holomorphic automorphic forms on unitary groups that are related to $\Pi$ can all be expressed as products of these $P_i(\Pi,\imath)$.  The Tate Conjecture predicts that all the {\bf unknown unknowns} have the same property.

\subsubsection*{Acknowledgements}\small
We thank Sug Woo Shin and Dipendra Prasad for several very useful conversations. We are also grateful to Nicola Antonio Porpora. HG also thanks the late Ferdinand Johannes G\"odde $\dagger$ (called ``Jan Loh'') for memorable discussions (``zero times zero'') in Bonn.  Finally, we are especially grateful to A. Raghuram, for his collaboration on the paper \cite{GHLR}, without which this present work would not have been possible, and for numerous conversations during the course of this work.
\normalsize

\numberwithin{equation}{section}
\section{Preliminaries} 

\subsection{Number fields and associate characters} \label{sect:fields}

We let $\Qbar$ be the algebraic closure of $\Q$ in $\C$. All number fields are considered as subfields of $\Qbar$.
For $k$ a number field, we let $J_k$ be its set of complex field-embeddings $\imath: k\hra\C$. We will write $S_\infty(k)$ for its set of archimedean places, $\O_k$ for its ring of integers, $\A_k$ for its ring of adeles, and use $k^{Gal}$ for a fixed choice of a Galois closure of $k/\Q$ in $\overline\Q\subset\C$. If $\pi$ is an abstract representation of a non-archimedean group, we will write $\Q(\pi)$ for the field of rationality of $\pi$, as defined in \cite{waldsp}, I.1. In this paper, every rationality field will turn out to be a number field.\\\\
Throughout our paper, $\cm$ will be reserved in order to denote a CM-field of dimension $2d = \dim_\Q \cm$. The set of archimedean places of $F$ is abbreviated $S_\infty=S_\infty(\cm)$. We will chose a section $S_\infty \ra J_F$ and may hence identify a place $v\in S_\infty$ with an ordered pair of conjugate complex embeddings $(\imath_v,\bar\imath_v)$ of $\cm$, where we will drop the subscript ``$v$'' if it is clear from the context. This order in turn fixes a choice of a CM-type $\Sigma:=\{\imath_v : v\in S_\infty\}$. The maximal totally real subfield of $\cm$ is denoted $\tr$. Its set of archimedean places will be identified with $S_\infty$, identifying a place $v$ with its first component embedding $\imath_v\in\Sigma$ and we let Gal$(\cm/\tr)=\{1,c\}$. \\\\
We extend the quadratic Hecke character $\varepsilon=\varepsilon_{\cm/\tr}: (\tr)^\times\backslash \A^\times_{\tr}\ra\C^\times$, associated to $\cm/\tr$ via class field theory, to a conjugate self-dual Hecke character $\eta: \cm^\times\backslash\A_{\cm}^\times\ra\C^\times$. At $v\in S_\infty$, $z\in F_v\cong\C$, we have $\eta_v(z)=z^t \bar z^{-t}$, where $t\in\tfrac12+\Z$. For the scope of this paper, we may assume without loss of generality that $t=\tfrac12$, \cite[\S 6.9.2]{bel-chen}. We define $\psi:=\eta\|\cdot\|^{1/2}$, which is an algebraic Hecke character;  here $\|\cdot\|$ is the normalized absolute value.\\\\
If $\chi$ is a Hecke character of $\cm$, we denote by $\widecheck{\chi}$ its conjugate inverse $(\chi^{c})^{-1}$.

\subsection{Algebraic groups and real Lie groups} \label{sect:alggrp}
 Let $(V_n,\<\cdot,\cdot\>)$ be an $n$-dimensional non-degenerate $c$-hermitian space over $\cm$, $n\geq 1$, we denote the corresponding unitary group over $F^+$ by $H:=H_n:=U(V_n)$. For each $v \in S_\infty$ we let $(r_v,s_v)$ denote the signature of the hermitian form induced by $\<\cdot,\cdot\>$ on the complex vector space $V_v := V\otimes_{F,\imath_v} \CC.$    \\\\ 
Whenever one has fixed an embedding $V_k\subseteq V_n$, we may view the attached unitary group $U(V_k)$ as a natural $\tr$-subgroup of $U(V_n)$. If $n=1$, the algebraic group $U(V_1)$ is isomorphic to the kernel of the norm map $N_{F/F^+}: R_{F/F^+}((\Gm)_F) \to (\Gm)_{F^+}$, where $R_{F/F^+}$ stands for the Weil-restriction of scalars from $F/F^+$, and is thus independent of $V_1$.\\\\ 
Let $\sigma\in{\rm Aut}(\C)$ and let $V_n$ be as above. Then there is a unique $c$-Hermitian space $^\sigma V_n$ over $F$, whose local invariants at the non-archimedean places of $F$ are the same as of $V_n$ and whose signatures satisfy $({}^\sigma r_v, {}^\sigma s_v)=(r_{\sigma^{-1}\circ \imath_v},s_{\sigma^{-1}\circ \imath_v})$ at all $v\in S_\infty$, cf.\ \cite{Lan}. We let ${}^\sigma H:= U({}^\sigma V_n)$ be the attached unitary group over $F^+$. By definition, ${}^\sigma H(\A_{f})\cong H(\A_{f})$ and ${}^\sigma H_\infty\cong \prod_{v\in S_\infty} H(F_{\sigma^{-1}\circ v})$.
\\\\
If $G$ is any reductive algebraic group over a number field $k$, we write $Z_G/k$ for its center, $G_\infty:=R_{k/\Q}(G)(\R)$ for the real Lie group of $\R$-points of the Weil-restriction of scalars from $k/\Q$ and denote by $K_{G,\infty}\subseteq G_\infty$ the product of $(Z_G)_\infty$ and a fixed choice of a maximal compact subgroup of $G_\infty$. Hence, we have $K_{G_{n,\infty}}\cong \prod_{v\in S_\infty}K_{G_n,v}$, each factor being isomorphic to $K_{G_n,v}\cong \R_+ U(n)$; $K_{H,\infty}\cong \prod_{v\in S_\infty}K_{H,v}$, with $K_{H,v}\cong U(r_v)\times U(s_v)$; and $K_{{}^\sigma H,\infty}\cong \prod_{v\in S_\infty}K_{H, \sigma^{-1}\circ v}$. Here, for any $m$, we denote by $U(m)$ the compact real unitary group of rank $m$.\\\\
Lower case gothic letters denote the Lie algebra of the corresponding real Lie group (e.g., $\g_{n,v}=Lie(G_n(\cm_v))$, $\k_{H,v}=Lie(K_{H,v})$, $\h_v=Lie(H(\tr_v))$, etc. ...). 

\subsection{Highest weight modules and cohomological automorphic representations}

\subsubsection{Finite-dimensional representations}  \label{sect:finitereps}

We let $\EE_\mu$ be an irreducible finite-dimensional representation of the real Lie group $G_{n,\infty}$ on a complex vector-space, given by its highest weight $\mu=(\mu_v)_{v\in S_\infty}$. 
Throughout this paper such a representation will be assumed to be algebraic: In terms of the standard choice of a maximal torus and a basis of its complexified Lie algebra, consisting of the functionals which extract the diagonal entries, this means that the highest weight of $\EE_\mu$ has integer coordinates, $\mu_v=(\mu_{\imath_v},\mu_{\bar\imath_v})\in\Z^n\times \Z^n$ for all $v\in S_\infty$. 
We say that $\EE_\mu$ is {\it m-regular}, if $\mu_{\imath_v,i}-\mu_{\imath_v,i+1}\geq m$ and $\mu_{\bar\imath_v,i}-\mu_{\bar\imath_v,i+1}\geq m$ for all $v\in S_\infty$ and $1\leq i\leq n-1$. 
Hence, $\mu$ is regular in the usual sense (i.e., inside the open positive Weyl chamber) if and only if it is $1$-regular.\\\\
Similarly, given a unitary group $H=U(V_n)$ we let $\cF_\lambda$ be an irreducible finite-dimensional representation of the real Lie group $H_\infty$ on a complex vector-space, given by its highest weight $\lambda=(\lambda_v)_{v\in S_\infty}$, $\lambda_v\in\Z^n$. Any such $\lambda$ may also be interpreted as the highest weight of an irreducible representation of $K_{H,\infty}$. In general, we will denote by $\Lambda=(\Lambda_v)_{v\in S_\infty}$ a highest weight for $K_{H,\infty}$ and we will write $\W_\Lambda$ for the corresponding irreducible representation.

\subsubsection{Cohomological representations} \label{sect:coh}

A representation $\Pi_\infty$ of $G_{n,\infty}$ is said to be {\it cohomological} if there is a highest weight module $\EE_\mu$ as above such that $H^*(\g_{n,\infty},K_{G_n,\infty},\Pi_\infty\otimes \EE_\mu)\neq 0$. In this case, $\EE_\mu$ is uniquely determined by this property and we say $\Pi_{\infty}$ is $m$-regular if $\EE_\mu$ is. \\\\
Analogously, a representation $\pi_\infty$ of $H_\infty$ is said to be {\it cohomological} if there is a highest weight module $\cF_\lambda$ as above such that $H^*(\h_\infty,K_{H,\infty},\pi_\infty\otimes \cF_\lambda)$ is non-zero. See \cite{bowa},\S I, for details.\\\\
It can be shown that an irreducible unitary generic representation $\Pi_\infty$ of $G_{n,\infty}$ is cohomological with respect to $\EE_\mu$ if and only if at each $v\in S_\infty$ it is of the form
\begin{equation}\label{eq:ind}
\Pi_v\cong {\rm Ind}_{B(\C)}^{G(\C)}[z_1^{a_{v,1}}\bar z^{-a_{v,1}}_1\otimes ...\otimes z_n^{a_{v,n}}\bar z^{-a_{v,n}}_n],
\end{equation}
where
\begin{equation}\label{hw and it}
a_{v,j}:=a(\mu_{\imath_v},j):=-\mu_{\imath_v,n-j+1}+\tfrac{n+1}{2}-j
\end{equation}
and induction from the standard Borel subgroup $B=TN$ is unitary, cf.\ \cite[Theorem 6.1]{enright} (See also \cite[\S 5.5]{grob-ragh} for a detailed exposition). The set $\{z^{a_{v,i}} \bar{z}^{-a_{v,i}}\}_{1\leq i\leq n}$ is called the {\it infinity type} of $\Pi_v$. For each $v$, the numbers $a_{v,i}\in \Z+\frac{n-1}{2}$ are all different and may be assumed to be in a strictly decreasing order, i.e. $a_{v,1}>a_{v,2}>\cdots>a_{v,n}$.\\\\
If $\pi_\infty$ is an irreducible tempered representation of $H_\infty$, which is cohomological with respect to $\cF_\lambda^{\sf v}$ (the presence of the contragredient will become clear in \S \ref{construction of motive}), then each of its archimedean component-representations $\pi_v$ of $H_v\cong U(r_v,s_v)$ is isomorphic to one of the $d_v:={n \choose r_v}$ inequivalent discrete series representations denoted $\pi_{\lambda,q}$, $0\leq q< d_v$, having infinitesimal character $\chi_{\lambda_v+\rho_v}$, \cite{vozu}. As it is well-known, \cite{bowa}, II Theorem 5.4, the cohomology of each $\pi_{\lambda,q}$ is centered in the middle-degree
$$H^p(\h_v,K_{H,v},\pi_{\lambda,q}\otimes\cF^{\sf v}_{\lambda_v})\cong\left\{\begin{array}{ll}
 \C & \textrm{if $p=r_vs_v$} \\
 0 & \textrm{else}
\end{array}
\right.$$
We thus obtain an $S_\infty$-tuple of Harish-Chandra parameters $(\HC_v)_{v \in S_\infty}$, and $\pi_\infty \cong \otimes_{v\in S_\infty} \pi_{\HC_v}$ where $\pi_{\HC_v}$ denotes the discrete series representation of $H_v$ with parameter $\HC_v$.

\subsubsection{Global base change and $L$-packets}\label{bc}

Let $\pi$ be a cohomological square-integrable automorphic\footnote{As usual, we will for convenience not distinguish between a square-integrable automorphic representation, its smooth limit-Fr\'echet-space completion or its (non-smooth) Hilbert space completion in the $L^2$-spectrum, cf.\ \cite{grob_zun} and \cite{grob_book} for a detailed account. Moreover, unless otherwise stated, an automorphic representation is always assumed to be irreducible.} representation of $H(\A_{F^+})$. It was first  proved by Labesse \cite{lab} (see also \cite{harris-labesse, kim-krish04, kim-krish05, morel,shin}) that $\pi$ admits a base change\footnote{Referring to \cite{shin}, the very careful reader may want to assume in addition to our standing assumptions on the field $F$ that $F=\mathcal K F^+$, where $\mathcal K$ is an imaginary quadratic field. This assumption, however, is resolved in \cite{KMSW}.} $BC(\pi) = \Pi$ to $G_n(\A_F)$: The resulting representation $\Pi$ is an isobaric sum $\Pi=\Pi_1\boxplus...\boxplus\Pi_k$ of conjugate self-dual square-integrable automorphic representations $\Pi_i$ of $G_{n_i}(\A_F)$, uniquely determined by the following: for every non-archimedean place $v$ of $F^+$, which splits in $F$ and where $\pi_v$ is unramified, the Satake parameter of $\Pi_v$ is obtained from that of $\pi_v$ by the formula for local base change, see for example \cite{ming}. 

It is then easy to see that at such places $v$, the local base change $\Pi_v$ is tempered if and only if $\pi_v$ is. The assumption that $\pi_{\infty}$ is cohomological implies moreover that $\Pi_{\infty}$ is cohomological: This was proved in \cite{lab} \S 5.1 for discrete series representations $\pi_\infty$ but follows in complete generality recalling that $\Pi_{\infty}$ has regular dominant integral infinitesimal character and hence is necessarily cohomological  by \cite{salamanca}, Theorem 1.8. It is then a consequence of \cite{bowa}, III.3.3 and the results in \cite{clozel, HT, Shin} -- here in particular \cite{car}, Theorem 1.2 -- that, if all isobaric summands $\Pi_i$ of $\Pi=BC(\pi)$ are cuspidal, all of their local components $\Pi_{i,v}$ are tempered. Here we also used the well-known fact that as the $\Pi_i$ are unitary, $\Pi$ is fully induced from its isobaric summands.\\\\
We define the global $L$-packet $\prod(H,\Pi)$ attached to such a representation $\Pi$ to be the set of cohomological tempered square-integrable automorphic representations $\pi$ of $H(\A_{F^+})$ such that $BC(\pi) = \Pi$. 
This is consistent with the formalism in \cite{mok, KMSW}, in which (as in Arthur's earlier work \cite{arthur}) the representation $\Pi$ plays the role of the global Arthur-parameter for the square-integrable automorphic representation $\pi$ of $H(\A_{F^+})$. We recall that temperedness together with square-integrability imply that $\pi$ is necessarily cuspidal, \cite{clozel2}, Proposition 4.10, \cite{wallach}, Theorem 4.3. Moreover, for  each $\pi\in \prod(H,\Pi)$, $\pi_\infty$ is in the discrete series, cf.\ \cite{vozu}. See also \cite{clozelihes}, Lemma 3.8 and Lemma 3.9. 

\begin{rmk}\label{rmk:temp}
It should be noted that for any cohomological cuspidal automorphic representation $\pi$ of $H(\A_{F^+})$, such that $\Pi=BC(\pi)$ is an isobaric sum $\Pi=\Pi_1\boxplus...\boxplus\Pi_k$ of conjugate self-dual cuspidal automorphic representations, $\pi_v$ is tempered at every place $v$ of $F^+$, i.e., in $\prod(H,\Pi)$. Indeed, in order to see this, recall that $\Pi$ serves as a generic, elliptic global Arthur-parameter $\phi$  in the sense of \cite{KMSW}, \S 1.3.4 (Observe that as $\Pi$ is cohomological, the isobaric summands must be all different.) Its localization $\phi_v$ at any place $v$ of $F^+$ (cf.\ \cite{KMSW}, Proposition 1.3.3), is bounded, because so is the local Langlands-parameter attached to the tempered representation $\Pi_{v}$ by the LLC, \cite{HT, henniart}. Hence, item (5) of Theorem 1.6.1 of \cite{KMSW} implies that each square-integrable automorphic representation $\pi$ of $H(\A_{F^+})$ attached to $\phi$ by \cite{KMSW}, Theorem 5.0.5 is tempered at all places. In particular, so is $\pi$.
\end{rmk}

\subsubsection{$\sigma$-twisted representations}
Let $\sigma\in{\rm Aut}(\C)$ and let $\Pi$ be a cohomological cuspidal automorphic representation of $G_n(\A_F)$. Then it is well-known that there exists a unique cohomological cuspidal automorphic representation ${}^\sigma\Pi$ of $G_n(\A_F)$, with the property that $({}^\sigma\Pi)_f\cong {}^\sigma(\Pi_f):=\Pi_f\otimes_{\sigma^{-1}}\C$, cf.\ \cite{clozel}, Theorem 3.13. Likewise, if $\pi$ is a cohomological cuspidal automorphic representation of $H(\A_{F^+})$, then there is a square-integrable automorphic representation ${}^\sigma\!\pi$ of ${}^\sigma H(\A_{F^+})$, such that $({}^\sigma\!\pi)_f\cong {}^\sigma(\pi_f):=\pi_f\otimes_{\sigma^{-1}}\C$: Recalling, \cite{gro-seb}, Theorem A.1 and \cite{milne-suh}, Theorem 1.3, this can be argued as in the second paragraph of \cite{BHR}, p.\ 665. In Lemma \ref{lem:unique} below we will provide conditions under which ${}^\sigma\!\pi$ is cuspidal and unique.

\subsection{Critical automorphic $L$-values and relations of rationality} \label{Galois equivariance}
 
 \subsubsection{Critical points of Rankin--Selberg $L$-functions}\label{sect:critRS}
Let $\Pi=\Pi_n\otimes\Pi_{n'}$ be the tensor product of two automorphic representations of $\GL_{n}(\A_F)\times\GL_{n'}(\A_F)$. We recall that a complex number $s_0\in \tfrac{n-{n'}}{2}+\Z$ is called {\it critical} for $L(s,\Pi_{n}\times\Pi_{n'})$ if both $L(s,\Pi_{{n},\infty}\times\Pi_{{n'},\infty})$ and $L(1-s,\Pi_{{n},\infty}^{\sf v}\times\Pi_{{n'},\infty}^{\sf v})$ are holomorphic at $s=s_0$. In particular, this defines the notion of critical points for standard $L$-functions $L(s,\Pi)$ and hence Hecke $L$-functions $L(s,\chi)$.\\\\
Let now $\Pi$ (resp. $\Pi'$) be a generic cohomological conjugate self-dual automorphic representation of $G_{n}(\Acm)$ (resp.\ $G_{n'}(\Acm)$) with infinity type $\{z^{a_{v,i}} \bar{z}^{-a_{v,i}}\}_{1\leq i\leq n}$ (resp.\ $\{z^{b_{v,j}} \bar{z}^{-b_{v,j}}\}_{1\leq j\leq n'}$) at $v\in S_\infty$. Then, the $L$-function $L(s,\Pi\times \Pi')$ has critical points if and only if $a_{v,i}+b_{v,j}\neq 0$ for all $v$, $i$ and $j$, cf.\ \S5.2 of \cite{jie-thesis}. In this case, the set of critical points of $L(s,\Pi\times\Pi')$ can be described explicitly as the set of numbers $s_0\in \tfrac{n-{n'}}{2}+\Z$ which satisfy
\begin{equation}
-\min |a_{v,i}+b_{v,j}| <s_{0}\leq \min |a_{v,i}+b_{v,j}|,
\end{equation}
the minimum being taken over all $1\leq i\leq n, 1\leq j\leq n'$, and $v\in S_{\infty}$. In particular, if $n\nequiv n' \mod 2$ then $s_0=\tfrac{1}{2}$ is always among these numbers. 

 \subsubsection{Relations of rationality and Galois equivariance}
 
 \begin{defn}[i]\label{definition algebraic relation}
Let  $E,L\subset\C$ be subfields and let $x, y\in E\otimes_\Q\C$. We write 
$$x\sim_{E\otimes_{\Q} L} y,$$ if either $y=0$, or, if $y$ is invertible and there is an $\ell\in E\otimes_\Q L$ such that $x=\ell y$ (multiplication being in terms of $\Q$-algebras). If the field $L$ equals $\Q$, it will be omitted in notation. \\
(ii) Let $E,L\subset\C$ again be subfields. Let $\underline x=\{x(\sigma)\}_{\sigma\in {\rm Aut}(\C)}$ and $\underline y=\{y(\sigma)\}_{\sigma\in {\rm Aut}(\C)}$ be two families of complex numbers. We write
$$\underline x\sim_{E} \underline y$$ 
and say that this relation {\it is equivariant under} Aut$(\C/L)$, if either $y(\sigma)=0$ for all $\sigma\in {\rm Aut}(\C)$, or if $y(\sigma)$ is invertible for all $\sigma\in {\rm Aut}(\C)$ and the following two conditions are verified:
\begin{enumerate}
\item $\cfrac{x(\sigma)}{y(\sigma)}\in \sigma(E)$ for all $\sigma$.
\item $\varrho\left(\cfrac{x(\sigma)}{y(\sigma)}\right)=\cfrac{x(\varrho\sigma)}{y(\varrho\sigma)} $ for all $\varrho\in {\rm Aut}(\C/L)$ and all $\sigma\in {\rm Aut}(\C)$.
\end{enumerate}
\end{defn}
Obviously, one may replace the first condition by requiring it only for all $\varrho$ running through representatives of ${\rm Aut}(\C)/{\rm Aut}(\C/L)$. In particular, if $L=\Q$, one only needs to verify it for the identity $id\in {\rm Aut}(\C)$. If moreover $E$ and $L$ are  number fields, one can define analogous relations for ${\rm Gal}(\Qbar/\QQ)$-families by replacing ${\rm Aut}(\C)$ by ${\rm Gal}(\Qbar/\QQ)$ and ${\rm Aut}(\C/L)$ by ${\rm Gal}(\Qbar/L)$. Note that a ${\rm Gal}(\Qbar/\QQ)$-family can be lifted to an ${\rm Aut}(\C)$-family via the natural projection ${\rm Aut}(\C)\rightarrow {\rm Gal}(\Qbar/\QQ)$, and two ${\rm Gal}(\Qbar/\QQ)$-families are equivalent if and only if their liftings are equivalent.

\begin{rem}[Aut$(\C)$-families vs. $\C^{|J_{E}|}$-tuples]\label{rem:Ealg}
Let $\underline x=\{x(\sigma)\}_{\sigma\in {\rm Aut}(\C)}$ and $\underline y=\{y(\sigma)\}_{\sigma\in {\rm Aut}(\C)}$ be two ${\rm Aut}(\C)$-families and assume we are given two number fields $E,L\subset \C$. If the individual numbers $x(\sigma)$, $y(\sigma)$ only depend on the restriction of $\sigma$ to $E$, then we may identify $\underline x$ and $\underline y$ canonically with elements $x,y\in\C^{|J_{E}|}\cong E\otimes_\Q\C$. The assertion that $\underline x\sim_{E} \underline y$, equivariant under Aut$(\C/L)$ implies that $x\sim_{E\otimes_{\Q} L} y$. 

Conversely, any element $x\in E\otimes_\Q\C \cong \C^{|J_{E}|}$ can be extended to a ${\rm Aut}(\C)$-family $\underline x=\{x(\sigma)\}_{\sigma\in {\rm Aut}(\C)}$, putting $x(\sigma):=x_{\sigma|_E}$. If we assume moreover that $E$ contains $L^{Gal}$, then for $x,y\in  E\otimes_\Q\C \cong \C^{|J_{E}|}$, the assertion $x\sim_{E\otimes_{\Q} L} y$ implies that $\underline x\sim_{E} \underline y$, equivariant under Aut$(\C/L)$.  We will be using this principle repeatedly in the precise form given in \cite[Lemma 1.34]{grob_lin}, in order to 
eliminate ambiguities introduced by auxiliary data.

In this paper, it will be convenient to have both points of view at hand. In fact, we prove assertions of the second type, which is generally a little bit stronger than the first one. But as we are always in the situation that $E$ contains $L^{Gal}$, the two assertions are equivalent and we will jump between them without further mention.  
\end{rem}

\subsection{Interlude: A brief review of motives and Deligne's conjecture}\label{general Deligne}

\subsubsection{Motives, periods over $\Q$ and Deligne's conjecture} We now quickly recall Deligne's conjecture about motivic $L$-functions, in order to put our main results into a precisely formulated framework and to fix notation. We follow Deligne, \cite{deligne}, \S 0.12, in adopting the following (common) pragmatic point of view through realizations:

\begin{defn}\label{definitionmotive}
A {\it motive} $M$ over a number field $k$ with {\it coefficients} in a number field $E(M)$ is a tuple
$$M=(M_{B,\imath}, M_{dR}, M_{\textrm{{\it \'et}}}; F_{B,\imath}, I_{\infty,\imath}, I_{\textrm{{\it \'et}},\imath}),$$
where $\imath\in J_k$ runs through the embeddings $k\hra\C$ and such that there exists an $n\geq 1$, where
\begin{enumerate}

\item[(B)] $M_{B,\imath}$ is an $n$-dimensional $E(M)$-vector space, together with a Hodge-bigraduation 
$$M_{B,\imath}\otimes_\Q \C = \bigoplus_{p,q} M^{p,q}_{B,\imath}$$
as a module over $E(M)\otimes_\Q\C$. 

\item[(dR)] $M_{dR}$ is a free $E(M)\otimes_\Q k$-module of rank $n$, equipped with a decreasing filtration $\{F^i_{dR}(M)\}_{i\in\Z}$ of $E(M)\otimes_\Q k$-submodules.

\item[(\'et)] $M_{\textrm{{\it \'et}}} = \{M_\ell\}_\ell$ is a strictly compatible system, cf.\ \cite{serrec} p.\ 11, of $\ell$-adic Gal$(\overline k/k)$-representations 
$$\rho_{M,\ell}: {Gal}(\overline k/k) \ra \GL(M_\ell)$$
on $n$-dimensional $E(M)_\ell$-vector spaces $M_\ell$, $\ell$ running through the set of finite places of $E(M)$,
\end{enumerate}
to be called ``realizations of $M$'', together with
\begin{enumerate}
\item[(i)] an $E(M)$-linear isomorphism 
$$F_{B,\imath}: M_{B,\imath} \ira M_{B,\overline\imath},$$
which satisfies $F_{B,\imath}^{-1} = F_{B,\overline\imath}$ and commutes with complex conjugation on the Hodge-bigraduation from (B), i.e., $\overline{F_{B,\imath}(M^{p,q}_{B,\imath})} \subseteq M^{p,q}_{B,\overline\imath}$,
\item[(ii)] an isomorphism of $E(M)\otimes_\Q \C$-modules
$$I_{\infty,\imath}: M_{B,\imath}\otimes_\Q\C \ira M_{dR}\otimes_{k,\imath}\C,$$
compatible with the Hodge-bigraduation from (B) and the decreasing filtration from (dR) above, i.e., $I_{\infty,\imath}(\bigoplus_{p\geq i} M^{p,q}_{B,\imath})=F^i_{dR}(M)\otimes_{k,\imath}\C$, and also compatible with $F_{B,\imath}$ and complex conjugation, i.e., $\overline{I_{\infty,\imath}} = I_{\infty,\overline\imath}\circ\overline{F_{B,\imath}}$, and  
\item[(iii)] a family $I_{\textrm{{\it \'et}},\imath} = \{I_{\imath,\ell}\}_\ell$ of isomorphisms of $E(M)_\ell$-vector spaces
$$I_{\imath,\ell}: M_{B,\imath}\otimes_{E(M)} E(M)_\ell \ira M_\ell,$$
$\ell$ running through the set of finite places of $E(M)$, where, if $\imath\in J_k$ is real, then $I_{\imath,\ell}\circ F_{B,\imath} = \rho_{M,\ell}(\gamma_\imath)\cdot I_{\imath,\ell}$, where $\gamma_\imath$ denotes complex conjugation of $\C$ attached to any extension to $\overline k$ of the embedding $\imath: k\hra \C$.
\end{enumerate}
to be called ``comparison isomorphisms''. The common rank $n$ of each realization as a free module is called the {\it rank} of $M$. If $n\geq 1$ and if there is an integer $w$ sucht that $M^{p,q}_{B,\imath}=\{0\}$ whenever $p+q\neq w$, then $M$ is called {\it pure of weight} $w$.
\end{defn}

The \'etale realization allows one to define the $E(M)\otimes_\Q\C$-valued $L$-{\it function} $L(s,M)$ of $M$ as the usual Euler product over the prime ideals $\p\lhd\O_k$,
$$L(s,M):=\left(\prod_{\p} L_\p(s,M)^\jmath\right)_{\jmath\in J_{E(M)}},$$
where $L_\p(s,M):=\det(id - N(\p)^{-s}\cdot \rho_{M,\ell}(Fr^{-1}_\p)| M^{I_\p}_{\ell})^{-1}$, and $Fr_\p$ denotes the geometric Frobenius locally at $\p$ (modulo conjugation) and $I_\p$ is the inertia subgroup in the decomposition group of an(y) extension of $\p$ to $\overline k$. Consequently, viewing $L_\p(s,M)$ as a rational function in the variable $X=N(\p)^{-s}$, the action of $\jmath\in J_{E(M)}$ on $L_\p(s,M)$ is defined by application to its coefficients: Here, we have to adopt the usual hypothesis, cf.\ \cite{deligne}, \S 1.2.1 \& \S 2.2, that at the finitely many ideals $\p$, where $\rho_{M,\ell}$ ramifies, the coefficients of $L_\p(s,M)$, viewed as a rational function in this way, belong to $E(M)$ and that they are independent of $\ell$ not dividing $N(\p)$, in order to obtain a well-defined element of $E(M)\otimes_\Q\C\cong \prod_{\jmath}\C$ (i.e., to make sense of the action of $\jmath$). It is well-known that $L(s,M)$ is absolutely convergent for $Re(s)\gg 0$ and it is tacitly assumed that $L(s,M)$ admits a meromorphic continuation to all $s\in\C$ as well as the usual functional equation with respect to the dual motive $M^{\sf v}$ (whose system of $\ell$-adic representations is contragredient to that of $M$), cf.\ \cite{deligne}, \S 2.2. An integer $m$ is then called {\it critical} for $L(s,M)$, if the archimedean $L$-functions on both sides of the functional equation are holomorphic at $s=m$. We refer to \cite{deligne}, \S 5.2 for the construction of the archimedean $L$-functions attached to $M$ and its dual. \\\\
Let now be $M$ a pure motive of weight $w$. By considering the motive $R_{k/\Q}(M)$, which is obtained from $M$ by applying restriction of scalars (i.e., whose system of $\ell$-adic representations is obtained by inducing the one attached to $M$ from Gal$(\overline k/k)$ to Gal$(\overline \Q/\Q)$) we may always reduce ourselves to the case, where $M$ is defined over $\Q$, which is the framework in which Deligne's conjecture is stated. As we are then left with only one embedding $\imath=id$, we will drag it along in order to lighten the burden of notation. \\\\
So, let $F_\infty=F_{B,id}: M_{B} \ira M_{B}$ be the only infinite Frobenius. If $w=2p$ is even, we suppose that $F_\infty$ acts by multiplication by $\pm 1$ on $M^{p,p}_B$. We then denote by $n^\pm=n^\pm(M)$ the dimension of the $+1$- (resp.\ $-1$-eigenspace) $M^\pm_B$ of $M_B$ of the involution $F_\infty$. Let $F_{dR}^\pm$ be $E(M)$-subspaces of $M_{dR}$, given by the filtration $\{F^i_{dR}(M)\}_{i\in\Z}$, such that the rank of $M^\pm_{dR}:=(M_{dR}/F^\mp_{dR})$ equals $n^\pm$ and such that $I_\infty$ induces isomorphisms of $E(M)\otimes_\Q\C$-modules
$$I^\pm_\infty: M^\pm_B\otimes_\Q\C \ira M^\pm_{dR}\otimes_\Q\C.$$
Following Delgine, we define two {\it periods}
$$c^\pm(M):=(\det(I^\pm_\infty)_\jmath)_{\jmath\in J_{E(M)}}\in (E(M)\otimes_\Q\C)^\times,$$
and
$$\delta(M):= (\det(I_\infty)_\jmath)_{\jmath\in J_{E(M)}}\in (E(M)\otimes_\Q\C)^\times.$$
Here, each determinant is computed with respect to a fixed choice of $E(M)$-rational bases of source and target spaces. Up to multiplication by an invertible element in the $\Q$-algebra $E(M)$, both periods hence depend only on $M$. 

\begin{conj}[Deligne, \cite{deligne}, Conjecture 2.8]\label{conj:Deligne}
Let $M$ be a pure motive of weight $w$ over $\Q$ and let $m$ be a critical point for $L(s,M)$. Then 
$$L^{S}(m,M)\sim_{E(M)} (2\pi i)^{n^{(-1)^m}\cdot m} \ c^{(-1)^m}(M)$$
\end{conj}

\subsubsection{Factorizing periods} 
Switching back to our general number field $k$, we choose and fix a section $S_\infty(k)\ra J_k$ and let $\Sigma_k$ be its image in $J_k$. If $w=2p$ is even, we assume, similar to the case $k=\Q$, that $R_{k/\Q}(\bigoplus_{\imath\in\Sigma_k} F_{B,\imath})$ acts by a scalar on $R_{k/\Q}(M)^{p,p}_B$. For $\imath\in \Sigma_k$ complex, this implies that $M^{p,p}_{B,\imath}=\{0\}$. Next, one may analogously define $\pm 1$--eigenspaces of $F_{B,\imath}$, which now obviously have to depend of the nature of the embedding $\imath\in J_k$: If $\imath$ is real, then our definition of $M_{B,\imath}^\pm$ is verbatim the one of the case $k=\Q$ from above, whereas if $\imath$ is complex, then we obtain eigenspaces $(M_{B,\imath}\oplus M_{B,\overline\imath})^\pm$ of the direct sum $M_{B,\imath}\oplus M_{B,\overline\imath}$. We also may analogously  define spaces $F^\pm_{dR}$ attached to the Hodge-filtration $\{F^i_{dR}(M)\}_{i\in\Z}$, cf.\ \cite{yoshida}, pp.\ 149--150, and we set $M^\pm_{dR}:=(M_{dR}/F^\mp_{dR})$. For $\imath\in\Sigma_k$, the maps $I_{\infty,\imath}$ induce canonical isomorphisms of $E(M)\otimes_\Q\C$-modules
$$I^\pm_{\infty,\imath}:M^\pm_{B,\imath}\otimes_\Q\C \ira M^\pm_{dR}\otimes_{k,\imath}\C,$$
if $\imath$ is real and 
$$I^\pm_{\infty,\imath}:(M_{B,\imath}\oplus M_{B,\overline\imath})^\pm\otimes_\Q\C \ira (M^\pm_{dR}\otimes_{k,\imath}\C)\oplus(M^\pm_{dR}\otimes_{k,\overline\imath}\C)$$
if $\imath$ is complex. 
We define
$$c^\pm(M,\imath):=(\det(I^\pm_{\infty,\imath})_\jmath)_{\jmath\in J_{E(M)}}\in (E(M)\otimes_\Q\C)^\times$$
and
$$\delta(M,\imath):= (\det(I_{\infty,\imath})_\jmath)_{\jmath\in J_{E(M)}} \in (E(M)\otimes_\Q\C)^\times.$$
Up to multiplication by an invertible element in $E(M)\otimes_\Q\imath(k)$, they only depend on $M$. Finally, let $n^\pm$ be the rank of the free $E(M)\otimes_\Q k$-module $M^\pm_{dR}$, if $k$ has a real place (respectively, if $k$ is totally imaginary, let $2n^\pm$ be the rank of the free $E(M)\otimes_\Q k$-module $M^\pm_{dR}\oplus \overline{M^\pm_{dR}}$, where $\overline{M^\pm_{dR}}$ is the $E(M)\otimes_\Q k$-module $M^\pm_{dR}$, but with complex conjugated scalar-mulitplication by $k$: $x\star v:=\overline x \cdot v$, $x\in k$, $v\in M^\pm_{dR}$.) Then, the two perspectives of Deligne's periods are linked by the following relations as elements of $\Q$-algebras:
$$c^\pm(R_{k/\Q}(M)) \sim_{E(M) K} D_k^{n^\pm/2}\prod_{\imath\in \Sigma_k} c^\pm(M,\imath)$$
$$\delta(R_{k/\Q}(M)) \sim_{E(M) K} D_k^{n/2} \prod_{\imath\in J_k} \delta(M,\imath),$$
where $K$ (resp.\ $D_k$) denotes the normal closure (resp.\ discriminant) of $k/\Q$, the latter identified with $1\otimes D_k$ in $E(M)\otimes_\Q \C$, cf.\ \cite{yoshida}, Proposition 2.2. We also refer to Proposition 2.11 of \cite{harris-lin} for a finer decomposition over $E(M)$.\\\\
Finally, we also recall the notion of regularity: To this end, we assume that we are given a motive $M$ with coefficients in $E(M)$. Since $E(M)\otimes_\Q\C\cong\C^{|J_{E(M)|}}$, for each $\imath\in J_k$,  there is a decomposition of $\C$-vector spaces
$$M^{p,q}_{B,\imath}=\bigoplus_{\jmath\in J_{E(M)}} M^{p,q}_{B,\imath}(\jmath).$$ 
We say that $M$ is {\it regular}, if $\dim M^{p,q}_{B,\imath}(\jmath)\leq 1$ for all $p,q\in\Z$, $\imath\in J_k$ and $\jmath\in J_{E(M)}$. For a fixed pair $(\imath,\jmath)$ as above, the set of pairs $(p,q)$, such that $M^{p,q}_{B,\imath}(\jmath)\neq 0$ is then called the {\it Hodge-type of $M$ at $(\imath,\jmath)$} and $(p,q)$ a {\it Hodge weight}. The Hodge-type is particularly useful, to give an explicit description of the critical points of $L(s,M)$. Indeed, if $M$ is regular and pure of even weight $w$, assume that $(\frac{w}{2},\frac{w}{2})$ is not a Hodge weight at any pair of embeddings. Then an integer $m$ is critical for $L(s,M)$ if and only if 
\begin{equation}
-\min_{(p,q)}\{ |p-\tfrac{w}{2}|\}+\tfrac{w}{2}< m\leq \min_{(p,q)}\{ |p-\tfrac{w}{2}|\}+\tfrac{w}{2}
\end{equation}
where $(p,q)$ runs over the Hodge weights for all pair $(\imath,\jmath)$.

\subsection{Motivic split indices}
Let $M$ and $M'$ be regular pure motives over $F$ with coefficients in a number field $E(M)=E(M')$ of weight $w$ and $w'$, respectively. We write $n$ for the rank of $M$ and $n'$ for the rank of $M'$. We write the Hodge-type of $M$ (resp.\ $M'$) at $(\imath,\jmath)$ as $(p_i,w-p_i)_{1\leq i\leq n}$, with $p_1>...>p_n$ (resp.\  $(q_j,w'-q_j)_{1\leq j\leq n'}$, with $q_1>...>q_{n'}$). Consider the tensor product $M\otimes M'$ (over $F$), whose system of $\ell$-adic representations is simply the system of tensor products $M_\ell\otimes M'_\ell$. We assume that $(M\otimes M')^{p,q}_{B,\imath}$ vanishes at $p=q=\tfrac{w+w'}{2}$, i.e., that $(\tfrac{w+w'}{2},\tfrac{w+w'}{2})$ is not a Hodge weight, i.e., $p_i+q_j\neq \tfrac{w+w'}{2}$ for all $i,j$. We put $p_{0}:=+\infty$ and $p_{n+1}:=-\infty$, and define:

$$sp(i,M;M',\imath,\jmath):=\#\{1\leq j\leq n'\mid p_{i}-\tfrac{w+w'}{2}>-q_{j}>p_{i+1}-\tfrac{w+w'}{2}\}.$$
We call $sp(i,M;M',\imath,\jmath)$ a {\it (motivic) split index}, reflecting the fact that the sequence of inequalities $-q_{n'}>...>-q_1$ splits into exactly $n+1$ parts, when merged with $p_1-\tfrac{w+w'}{2}>...>p_n-\tfrac{w+w'}{2}$, where the length of the $i$-th part in this splitting is $sp(i,M;M',\imath,\jmath)$. This gives rise to the following

\begin{defn}\label{split motivic}
For $0\leq i\leq n$ and $\imath\in\Sigma$, we define the {\it (motivic) split indices} (cf.\ \cite{harris-lin}, Definition 3.2)
$$sp(i,M;M',\imath):=(sp(i,M;M',\imath,\jmath))_{\jmath\in J_{E(M)}}\in\N^{J_{E(M)}},$$
and, mutatis mutandis, 
$$sp(j,M';M,\imath):=(sp(j,M';M,\imath,\jmath))_{\jmath\in J_{E(M')}}\in\N^{J_{E(M')}},$$
\end{defn}

\subsection{Motivic periods}\label{motivic periods}

Let $M$ be a regular pure motive over $\cm$ of rank $n$ and weight $w$ with coefficients in a number field $E\supset F^{Gal}$. For $1\leq i\leq n$ and $\imath \in \Sigma$, we have defined {\it motivic periods} $Q_{i}(M,\imath)$ in \cite{harrisANT} (see \cite{harris-lin}, Definition 3.1 for details). They are elements in $E\otimes_\Q \C$, well-defined up to multiplication by elements in $E\otimes_\Q \imath(\cm)$. If $M$ is moreover polarised, i.e., if $M^{{\sf v}}\cong M^{c}$, the period $Q_{i}(M,\imath)$ is equivalent to the inner product of a vector in $M_{B,\imath}$, the Betti realisation of $M$ at $\imath$, whose image via the comparison isomorphism is inside $i$-th bottom degree of the Hodge filtration for $M$. We have furthermore defined
\begin{equation}\label{eq:Qperiods}
Q^{(i)}(M,\imath):=Q_{0}(M,\imath)Q_{1}(M,\imath)\cdots Q_{i}(M,\imath), 
\end{equation}
where $Q_{0}(M,\imath):=\delta(M,\imath)(2\pi i)^{n(n-1)/2}$. Then, Deligne's periods can be interpreted interpreted in terms of the above motivic periods:

\begin{prop}(cf.\ \cite{harris-lin}, Proposition 2.11 and 3.13)\label{Deligne period motivic}
Let $M$ be a regular pure motive over $\cm$ of rank $n$ and weight $w$ with coefficients in a number field $E\supset F^{Gal}$ and let $M'$ be a regular pure motive over $\cm$ of rank $n'$ and weight $w'$ with coefficients in a number field $E'\supset F^{Gal}$. We assume that $(\tfrac{w+w'}{2},\tfrac{w+w'}{2})$ is not a Hodge weight for the motive $M\otimes M'$ with coefficients in $EE'$. Then, the Deligne periods satisfy
\begin{eqnarray}
&c^{\pm}(R_{\cm/\Q}(M\otimes M'))&\\ \nonumber
 &\sim_{EE'\otimes_\Q F^{Gal}}  (2\pi i)^{-\frac{nn'd(n+n'-2)}{2}} \prod\limits_{\imath\in \Sigma} [\prod\limits_{j=0}^{n}Q^{(j)}(M,\imath)^{sp(j,M;M',\imath)}\prod\limits_{k=0}^{n'}Q^{(k)}(M',\imath)^{sp(k,M';M,\imath)}].
\end{eqnarray} 
\end{prop}

\section{Translating Deligne's conjecture into an automorphic context}\label{auto Deligne}

\subsection{CM-periods and special values of Hecke characters}\label{CM-periods}
\subsubsection{Special Simura data and CM-periods}\label{sect:CM}
Let $(T,h)$ be a Shimura datum where $T$ is a torus defined over $\Q$ and $h:R_{\C/\R}(\mathbb{G}_{m,\C})\rightarrow T_{\R}$ a homomorphism satisfying the axioms defining a Shimura variety, cf.\ \cite{milne}, II. Such pair is called a {\it special} Shimura datum. Let $Sh(T,h)$ be the associated Shimura variety and let $E(T,h)$ be its reflex field.\\\\
For $\chi$ an algebraic Hecke character of $T(\A_\Q)$, we let $E_T(\chi)$ be the number field generated by the values of $\chi_f$, $E(T,h)$ and $F^{Gal}$, i.e., the composition of the rationality field $\Q(\chi_f)$ of $\chi_f$, and $E(T,h) F^{Gal}$. If it is clear from the context, we will also omit the subscript ``$T$''. We may define a non-zero complex number $p(\chi,(T,h))$, called {\it CM-period}, as in Sect.\ $1$ of \cite{Harris93} and the appendix of \cite{harrisappendix}, to which we refer for details: It is defined as the ratio between a certain deRham-rational vector and a certain Betti-rational vector inside the cohomology of the Shimura variety with coefficients in a local system. As such, it is well-defined modulo $E_T(\chi)^{\times}$. Recall the $\sigma$-twisted Shimura datum, $({}^\sigma T,{}^\sigma h)$, $\sigma\in {\rm Aut}(\C)$, cf.\ \cite{milne}, II.4, \cite{langlands}. By taking ${\rm Aut}(\C)$-conjugates of the aforementioned rational vectors, we can define the family $\{p({}^\sigma\chi,({}^\sigma T,{}^\sigma h))\}_{\sigma\in {\rm Aut}(\C)}$, such that if $\sigma$ fixes $E_{T}(\chi)$ then $p({}^\sigma\chi,({}^\sigma T,{}^\sigma h))=p(\chi,(T,h))$, i.e., in view of Remark \ref{rem:Ealg}, $\{p({}^\sigma\chi,({}^\sigma T,{}^\sigma h))\}_{\sigma\in {\rm Aut}(\C)}$ defines an element in $\C^{|J_{E_T(\chi)}|}\cong E_T(\chi)\otimes_\Q\C$. The following proposition holds ${\rm Aut}(\C)$-equivariantly as interpreted for the family $\{p({}^\sigma\chi,({}^\sigma T,{}^\sigma h))\}_{\sigma\in {\rm Aut}(\C)}$:\\

\begin{prop}\label{relation CM-period}
Let $T$ and $T'$ be two tori defined over $\Q$ both endowed with a special Shimura datum $(T,h)$ and $(T',h')$ and let $u:(T',h')\rightarrow (T,h)$ be a homomorphism between them. Let $\chi$ be an algebraic Hecke character of $T(\A_\Q)$ and put $\chi':=\chi\circ u$, which is an algebraic Hecke character of $T'(\A_\Q)$. Then we have:
\begin{equation}\nonumber
p(\chi,(T,h)) \sim_{E_T(\chi)} p(\chi',(T',h')).
\end{equation}
Interpreted as families, this relation is equivariant under the action of ${\rm Aut}(\C)$. 
\end{prop}
\begin{proof} 
This is due to the fact that both the Betti-structure and the deRham-structure commute with the pullback map on cohomology. We refer to \cite{Harris93}, in particular relation $(1.4.1)$ for details.
\end{proof}
If $\Psi$ a set of embeddings of $F$ into $\C$ such that $\Psi\cap \Psi^{c}=\emptyset$, one can define a special Shimura datum $(T_{F},h_{\Psi})$ where $T_{F}:=R_{F/\Q}(\mathbb{G}_{m})$ and $h_{\Psi}:R_{\C/\R}(\mathbb{G}_{m,\C}) \rightarrow T_{F,\R}$ is a homomorphism such that over $\imath\in J_{F}$, the Hodge structure induced by $h_{\Psi}$ is of type $(-1,0)$ if $\imath\in \Psi$, of type $(0,-1)$ if $\imath\in \Psi^{c}$, and of type $(0,0)$ otherwise. In this case, for $\chi$ an algebraic Hecke character of $F$, we write $p(\chi,\Psi)$ for $p(\chi, (T_{F},h_{\Psi}))$ and abbreviate $p(\chi, \imath):=p(\chi,\{\imath\})$. We also define the (finite) compositum of number fields $E_F(\chi):=\prod_\Psi E_{T_F}(\chi)$.

\begin{lem}\label{Lemma CM} 
Let $\imath\in \Sigma$ and let $\Psi$ and $\Psi'$ be disjoint sets of embeddings of $F$ into $\C$ such that $\Psi\cap \Psi^{c}=\emptyset= \Psi'\cap \Psi'^{c}$. Let $\chi$ and $\chi'$ be algebraic Hecke characters of $\GL_1(\A_F)$, 
and recall the algebraic Hecke character $\psi$ from Sect.\ \ref{sect:fields}. Then,
\begin{enumerate}[label=(\alph*)] 
\item $p(\chi, \Psi\sqcup \Psi') \sim_{E_F(\chi)} p(\chi,\Psi) \ p(\chi, \Psi')$
\item  $p(\chi\chi',\Psi) \sim_{E_F(\chi)E_F(\chi')} p(\chi,\Psi) \ p(\chi',\Psi) $
\item If $\chi$ is conjugate selfdual, then $p(\widecheck{\chi},\bar{\imath})\sim_{E_F(\chi)} p(\widecheck{\chi},\imath)^{-1}.$ 
\item $p(\widecheck{\psi},\bar{\imath})\sim_{E_F(\psi)} (2\pi i)p(\psi,\imath)^{-1}$.
\end{enumerate}
Interpreted as families, these relations are equivariant under the action of ${\rm Aut}(\C)$.
\end{lem}
\begin{proof}
The first two assertions are proved in \cite{grob_lin}, Proposition 4.4. For (c), observe that by Lemma $1.6$ of \cite{Harris93}, we have $p(\widecheck{\chi},\bar{\imath})\sim_{E_F(\chi)} p(\widecheck{\chi}^{c},\imath)$. Then Proposition $1.4$ of \cite{Harris93} and the fact that $\chi \chi^{c}$ is trivial imply
$$p(\widecheck{\chi}^{c},\imath)p(\widecheck{\chi},\imath)\sim_{E_F(\chi)}p(\widecheck{\chi \chi^{c}},\imath)\sim_{E_F(\chi)}1.$$ 
Similarly, for the last assertion, we have
$p(\widecheck{\psi},\bar{\imath})\sim_{E_F(\psi)} p(\widecheck{\psi}^{c},\imath)\sim_{E_F(\psi)} p(\widecheck{\psi}^{-1}\|\cdot \|^{-1},\imath)\sim_{E_F(\psi)}(2\pi i)p(\psi,\imath)^{-1}$ where the last step is due to the fact that $p(\|\cdot\|,\imath)\sim_{\Q}(2\pi i)^{-1}$ (cf. $1.10.9$ of \cite{H97}). 
\end{proof}



\subsection{Arithmetic automorphic periods}\label{arithmetic automorphic period}
\subsubsection{A theorem of factorization}

In this paper we focus on cohomological conjugate self-dual, cuspidal automorphic representations $\Pi$ of $\GL_{n}(\Acm)$, which satisfy the following assumption: 

\begin{hyp}\label{descent} 
For each $I=(I_{\imath})_{\imath\in \Sigma}\in \{0,1,\cdots,n\}^{|\Sigma|}$ there is a unitary group $H_I$ over $F^{+}$ as in \S \ref{sect:alggrp} of signature $(n-I_\imath,I_\imath)$ at $v=(\imath,\bar\imath)\in S_\infty$ such that the global $L$-packet $\prod(H_I,\Pi)$ is non-empty. Moreover, if $\pi\in \prod(H_I,\Pi)$, then the packet also contains all the representations $\tau_\infty\otimes\pi_f$, $\tau_\infty$ running through the discrete series representation of $H_{I,\infty}$ of the same infinitesimal character of $\pi_\infty$.
\end{hyp}

\begin{rmk}\label{rmk:des}
This hypothesis is always satisfied, if $n$ is odd. For $n$ even it is also known to hold, if $\Pi_\infty$ is cohomological with respect to a regular representation and $\Pi_v$ is square-integrable at a non-archimedean place $v$ of $F^+$, which is not split in $F$. Moreover, it is well-known that a cohomological conjugate self-dual, cuspidal automorphic representation of $\GL_{n}(\Acm)$ always descends to a cohomological cuspidal automorphic representation of the quasi-split unitary group $U^*_n$ of rank $n$ over $F^+$, cf.\ \cite{harris-labesse} and \cite{mok}, Corollary 2.5.9 (and the argument in \cite{grob_harris_lapid}, \S 6.1). In contrast to this positive result, there are also cohomological conjugate self-dual, cuspidal automorphic representations $\Pi$ of $\GL_{n}(\Acm)$, which do not satisfy Hypothesis \ref{descent}: As the simplest counterexample, take an everywhere unramified Hilbert modular cusp form for a real quadratic field $F^+$ not of CM-type. The quadratic base change of the corresponding automorphic representation to a CM-quadratic extension $F$ does not descend to a unitary group of signature $(1,1)$ at one archimedean place and $(2,0)$ at the other.
\end{rmk}
If $\Pi$ satisfies Hypothesis \ref{descent}, a family of \textit{arithmetic automorphic periods} $\{P^{({}^{\sigma}I)}({}^{\sigma}\Pi)\}_{\sigma\in {\rm Aut}(\C)}$ can then be defined as the Petersson inner products of an ${\rm Aut}(\C)$-equivariant family of arithmetic holomorphic automorphic forms as in (2.8.1) of \cite{H97}, or Definition 4.6.1 of \cite{jie-thesis}. The following result is proved in \cite{linfactorization}, Theorem 3.3.

\begin{thm}[Local arithmetic automorphic periods]\label{thm:artihap} 
Let $\Pi$ be a cohomological conjugate self-dual cuspidal automorphic representation of $\GL_{n}(\Acm)$, which satisfies Hypothesis \ref{descent}. We assume that  $\Pi$ is $5$-regular.  
Then there exists a number field $E(\Pi)\supseteq \Q(\Pi_f) F^{Gal}$ (see \S \ref{sect:EPi} below) and families of {\rm local} arithmetic automorphic periods $\{P^{(i)}({}^{\sigma}\Pi,\imath)\}_{\sigma\in {\rm Aut}(\C)}$, for $0\leq i \leq n$ and $\imath\in\Sigma$, which are unique up to multiplication by elements in $E(\Pi)^\times$ such that \begin{equation}\label{local end}
P^{(0)}(\Pi,\imath)\sim_{E(\Pi)} p(\widecheck{\xi}_{\Pi},\bar\imath) \quad\quad {\it and} \quad\quad P^{(n)}(\Pi,\imath)\sim_{E(\Pi)} p(\widecheck{\xi}_{\Pi},\imath),
\end{equation} 
where $\xi_{\Pi}$ denotes the central character of $\Pi$, and satisfy the relation
\begin{equation}\label{eq:splitting}
 P^{(I)}(\Pi) \sim_{E(\Pi)} \prod\limits_{\imath\in\Sigma}P^{(I_{\imath})}(\Pi,\imath).
 \end{equation}
 In particular, we have \begin{equation}\label{local end 2}
P^{(0)}(\Pi,\imath) P^{(n)}(\Pi,\imath)\sim_{E(\Pi)} 1.
\end{equation} 
Interpreted as families, all relations are equivariant under the action of ${\rm Aut}(\C/F^{Gal})$. 
\end{thm}



\subsubsection{The field $E(\Pi)$}\label{sect:EPi} 
Let $\pi\in \prod(H_I,\Pi)$. It is easy to see that the {\it field of rationality} $\QQ(\pi_f)$ of $\pi_f$, which, as we recall, is defined as the fixed field in $\C$ of the subgroup of $\sigma \in {\rm Aut}(\C)$ such that ${}^\sigma\pi_f \cong \pi_f$, coincides with $\QQ(\Pi_f)$ if the local $L$-packets that base change to $\Pi$ are singletons, and are simple finite extensions of $\QQ(\Pi_f)$ otherwise, see Lemma \ref{contained}.   However, because of the presence of non-trivial Brauer obstructions it is not always possible to realize $\pi_f$ over $\QQ(\pi_f)$.  By {\it field of definition} we mean a field over which $\pi_f$ has a model. The field $E(\Pi)$ in the statement of Theorem \ref{thm:artihap} may be taken to be the compositum of $F^{Gal}$ with fields of definition of the descents $\pi\in \prod(H_I,\Pi)$, $I=(I_{\imath})_{\imath\in \Sigma}\in \{0,1,\cdots,n\}^{|\Sigma|}$, as in Hypothesis \ref{descent}, cf.\ \cite{linfactorization}, Theorem 2.2.: It follows from \cite{gro-seb2}, Theorem A.2.4, that these (finitely many) fields of definition exist and are number fields. In fact, they can be taken to be finite abelian extension of the respective $\QQ(\pi_f)$ {\it From now on, $E(\Pi)$ will stand for (any fixed choice of) such a field.}

\subsection{Automorphic split indices}
Let $n$ and $n'$ be two integers. Let $\Pi$ (resp. $\Pi'$) be a cohomological conjugate self-dual cuspidal automorphic representation of $G_{n}(\Acm)$ (resp.\ $G_{n'}(\Acm)$) with infinity type $\{z^{a_{v,i}} \bar{z}^{-a_{v,i}}\}_{1\leq i\leq n}$ (resp.\ $\{z^{b_{v,j}} \bar{z}^{-b_{v,j}}\}_{1\leq j\leq n'}$) at $v\in S_\infty$.

\begin{defn}\label{split automorphic}
For $0\leq i\leq n$ and $\imath_v\in\Sigma$, we define the {\it automorphic split indices}, cf.\ \cite{jie-thesis,harris-lin},
$$sp(i,\Pi;\Pi',\imath_v):=\#\{1\leq j\leq n'\mid -a_{v,n+1-i}>b_{v,j}>-a_{v,n-i}\}$$
and
$$sp(i,\Pi;\Pi',\bar\imath_v):=\#\{1\leq j\leq n' \mid a_{v,i}>-b_{v,j}>a_{v,i+1}\}.$$ 
\noindent Here we put formally $a_{v,0}=+\infty$ and $a_{v,n+1}=-\infty$. It is easy to see that
\begin{equation}\label{sp relation}
sp(i,\Pi^{c};\Pi'^{c},\imath_v)=sp(i,\Pi;\Pi',\bar{\imath}_v)=sp(n-i,\Pi;\Pi',\imath_{v}).
\end{equation}
Similarly, for $0\leq j\leq n'$, we define $sp(j,\Pi';\Pi,\imath_v):=\#\{1\leq i\leq n\mid -b_{v,n'+1-j}>a_{v,i}>-b_{v,n'-j}\}$ and $sp(j,\Pi';\Pi,\bar\imath_v):=\#\{1\leq i\leq n\mid b_{v,j}>-a_{v,i}>b_{v,j+1}\}$.
\end{defn}

\subsection{Translating Deligne's conjecture for Rankin--Selberg $L$-functions}\label{tensorp}

We resume the notation and assumptions from the previous section and we suppose moreover that $a_{v,i}+b_{v,j}\neq 0$ for any $v\in S_\infty$, $1\leq i\leq n$ and $1\leq j\leq n'$. \\\\
Conjecturally, there are motives $M=M(\Pi)$ (resp. $M'=M(\Pi')$) over $F$ with coefficients in a finite extension $E$ of $\Q(\Pi_f)$ (resp. $E'$ of $\Q(\Pi'_f)$), satisfying $L(s,R_{\cm/\Q}(M\otimes M'))=L(s-\tfrac{n+n'-2}{2},\Pi_f\times \Pi'_f)$, which is a variant of \cite{clozel}, Conjecture 4.5. To make sense of this statement, the right hand side of the equation must first be interpreted as a function with values in $EE'\otimes_\Q\C\cong \C^{|J_{EE'}|}$: Arguing as in \cite{grob_harris_lapid}, \S 4.3, or in \cite{clozel}, Lemma 4.6, one shows that at $v\notin S_\infty$, the local $L$-factor $L(s-\frac{n+n'-2}{2},\Pi_v\times \Pi'_v)=P_v(q^{-s})^{-1}$ for a polynomial $P_v(X)\in EE'[X]$, satisfying $P(0)=1$, and one deduces that $L(s-\frac{n+n'-2}{2},{}^\sigma\Pi_v\times {}^\sigma\Pi'_v)={}^\sigma P_v(q^{-s})^{-1}$, where $\sigma$ acts on $P_v$ by application to its coefficients in $EE'$. In particular, for any finite set $S$ of places of $F$ containing $S_\infty$, the family $\{\prod_{v\notin S} L(s-\frac{n+n'-2}{2},{}^\sigma\Pi_v\times {}^\sigma\Pi'_v)\}_{\sigma\in {\rm Aut}(\C)}$ only depends on the restriction of the individual $\sigma$ to $EE'$, whence we may apply Remark \ref{rem:Ealg}, in order to view it as an element of $\C^{|J_{EE'}|}\cong EE'\otimes_\Q\C$. It is this way, in which we will interpret $L^S(s-\tfrac{n+n'-2}{2},\Pi\times \Pi')$ as a $|J_{EE'}|$-tuple.\\\\
In \S \ref{construction of motive} we shall indeed construct such motives $M=M(\Pi)$ and $M'=M(\Pi')$, in the sense of collections of cohomological realizations, attached to a large family of representations $\Pi$ and $\Pi'$,  (consistently with Clozel's conjecture mentioned above).  It will turn out that $M$ (resp. $M'$) is regular, pure of rank $n$ (resp. $n'$) and weight $w=n-1$ (resp. $w'=n'-1$) whose field of coefficients may be chosen to be a suitable finite extension $E$ of $E(\Pi)$, resp.\ $E'$ of $E(\Pi')$. Moreover, the above condition on the infinity type is equivalent to the condition that the $(\tfrac{w+w'}{2},\tfrac{w+w'}{2})$ Hodge component of $M\otimes M'$ is trivial. We can hence apply Proposition \ref{Deligne period motivic} and obtain a relation between Deligne's periods $c^\pm(R_{F/\Q}(M\otimes M'))$ and our motivic periods $Q^{(i)}(M,\imath)$ and $Q^{(j)}(M',\imath)$.\\\\
It is predicted by the Tate conjecture (see Conjecture 2.8.3 and Corollary 2.8.5 of \cite{H97} and Sect.\ $4.4$ of \cite{harris-lin}), that one has the fundamental {\it Tate-relation}: 
\begin{equation}\label{eq:fundrel}
P^{(i)}(\Pi,\imath) \sim_{E} Q^{(i)}(M(\Pi),\imath).
\end{equation}
Again, as for the comparison of motivic and automorphic $L$-functions above, the left hand side of this relation should be read as an element of $E\otimes_\Q\C$ as explained in Remark \ref{rem:Ealg}. Let us assume for a moment that \eqref{eq:fundrel} is valid (and that its left hand side is defined).\\\\ One easily checks that our so constructed motive $M$ has Hodge type at $\imath$ given by
$$(-a_{v,n+1-i}+w/2,a_{v,n+1-i}+w/2)_{1\leq i\leq n}, \text{ if } \imath = \imath_v;~~  (a_{v,i} + w/2, -a_{v,i} + w/2)_{1 \leq i \leq n}, \text{ if } \imath = \bar{\imath}_v.$$
Similarly, $M'$ has Hodge type at $\imath$ given by 
 $$(-b_{v,n'+1-j}+w'/2,b_{v,n'+1-j}+w'/2)_{1\leq j\leq n'} \text{ if } \imath = \imath_v; ~~(b_{v,j} + w'/2, -b_{v,j} + w'/2)_{1 \leq j \leq n'}, \text{ if } \imath = \bar{\imath}_v.$$ 
We now see immediately from Definition \ref{split motivic} and Definition \ref{split automorphic} that $sp(i,\Pi;\Pi',\imath)=sp(i,M;M',\imath)$. Whence, recalling that for any critical point $s_0\in \tfrac{n+n'}{2}+\Z$ of $L(s,\Pi\times \Pi')$, and any $v\notin S_\infty$, $L(s_0,\Pi_v\times \Pi'_v)$ is the inverse of a polynomial expression $P_v(q^{-s})\in EE'[q^{-s}]$ of an integral power $s_0-\tfrac{n+n'}{2}$ of $q$, and hence in $EE'$, and recollecting all of our previous observations, we finally deduce that Deligne's conjecture, Conjecture \ref{conj:Deligne}, for $R_{\cm/\Q}(M\otimes M')$ may be rewritten in purely automorphic terms as follows:

\begin{conj}\label{main conjecture}
Let $\Pi$ (resp. $\Pi'$) be a cohomological conjugate self-dual cuspidal automorphic representation of $G_{n}(\Acm)$ (resp. $G_{n'}(\Acm)$), which satisfies Hypothesis \ref{descent}. Let $s_{0}\in \Z+\tfrac{n+n'}{2}$ be a critical point of $L(s,\Pi\times \Pi')$, and let $S$ be a fixed finite set of places of $F$, containing $S_\infty$. Then, the arithmetic automorphic periods $P^{(I)}(\Pi)$ and $P^{(I)}(\Pi')$ admit a factorization as in \eqref{eq:splitting} and
\begin{equation}
L^S(s_{0},\Pi\otimes \Pi') \sim_{E(\Pi)E(\Pi')} (2\pi i)^{nn's_{0}} \prod\limits_{\imath \in \Sigma}[\prod\limits_{0\leq i\leq n}P^{(i)}(\Pi,\imath)^{sp(i,\Pi;\Pi',\imath)}\prod\limits_{0\leq j\leq n'}P^{(j)}(\Pi',\imath)^{sp(j,\Pi';\Pi,\imath)}].
\end{equation}
Interpreted as families, this relation is equivariant under action of ${\rm Aut}(\C/F^{Gal})$.
\end{conj}

\subsection{About the main goals of this paper and a remark on the strategy of proof}\label{sect:goals}

In \cite{GHLR} the authors, together with Raghuram, proved the automorphic version 
Theorem \ref{automorphic Deligne general intro} of Deligne's Conjecture 
for the tensor products of motives $M(\Pi)\otimes M(\Pi')$ attached (as in Theorem \ref{thm:clozel}) to a large family of cohomological conjugate self-dual cuspidal automorphic representations $\Pi$ and $\Pi'$ of $G_{n}(\Acm)$, resp.\ $G_{n'}(\Acm)$.  We used  
Theorem \ref{thm:artihap} to rewrite the initial automorphic formula, and thus to verify Conjecture \ref{main conjecture} for those $\Pi$ and $\Pi'$. 
\\\\ 
In view of Proposition \ref{Deligne period motivic} this will reduce a complete proof of Deligne's original conjecture, Conjecture \ref{conj:Deligne}, for the motives attached to such $\Pi$ and $\Pi'$ (and with coefficients in a number field containing $F^{Gal}$) to a proof of the Tate relation \eqref{eq:fundrel}.  
Our main result is a version of \eqref{eq:fundrel} by showing a refined decomposition of the local arithmetic automorphic periods $P^{(i)}(\Pi,\imath)$, which mirrors \eqref{eq:Qperiods}: Recall that the motivic periods on the right-hand-side of the Tate relation were defined as a product
$$Q^{(i)}(M(\Pi),\imath)=Q_{0}(M(\Pi),\imath)Q_{1}(M(\Pi),\imath)\cdots Q_{i}(M(\Pi),\imath).$$
For our main result we will define factors $P_i(\Pi,\imath)$ in \S\ref{sect:thm}, which are attached to a certain (canonical) descent $\pi(i)$ of $\Pi$ to a (non-canonical) unitary group and show that, up to a scalar contained in an extension $E\supset F^{Gal}$ explicitly attached to $\Pi$ and the $\pi(i)$'s, we have
\begin{equation}\label{eq:factor}
P^{(i)}(\Pi,\imath) \sim_{E} P_{0}(\Pi,\imath)P_{1}(\Pi,\imath)\cdots P_{i}(\Pi,\imath).
\end{equation}
The $P_i$ is essentially (but not quite) the {\it automorphic $Q$-period} of  $\pi(i)$ introduced in \S \ref{Qperiod}, and the Tate relation then comes down to a rather simple comparison, established in and recorded as Theorem \ref{main factorization}.\\\\

\section{Shimura varieties, coherent cohomology and a motive}

\subsection{Shimura varieties for unitary groups}\label{Sdata}
Let $V=V_n$ and $H=U(V)$ be as defined in \S \ref{sect:alggrp}. Let $\Ss = R_{\CC/\RR} \mathbb{G}_{m,\CC}$, so that $\Ss(\RR) = \CC^\times$, canonically.  In this paper we will use period invariants, attached to a Shimura datum $(H,Y_V)$, as in \cite[\S 2.2]{H21}. Explicitly, the base point $y_V  \in Y_V$ is given by 

\begin{equation}\label{yV} y_{V,v}(z) = \begin{pmatrix}  (z/\bar{z})I_{r_v} & 0 \\ 0 & I_{s_v} \end{pmatrix} \end{equation}
The following lemma is then obvious: We record it here in order to define parameters for automorphic vector bundles in the next sections.  

\begin{lem}\label{stab}  
Let $y \in Y_V$.   Its stabilizer $K_y=:K_{H,\infty}$ in $H_\infty$ is isomorphic to $\prod_{v \in S_\infty} U(r_v)\times U(s_v)$.
\end{lem}
Unlike the Shimura varieties attached to unitary similitude groups, the Shimura variety $Sh(H,Y_V)$ attached to $(U(V),Y_V)$ parametrizes Hodge structures of weight $0$ -- the homomorphisms $y \in Y_V$ are trivial on the subgroup $\RR^\times \subset \CC^\times$ -- and are thus of abelian type but not of Hodge type. The reflex field $E(H,Y_V)$ is the subfield of $F^{Gal}$ determined as the stabilizer of the cocharacter $\kappa_V$ with $v$-component  $\kappa_{V,v}(z) = \begin{pmatrix}  zI_{r_v} & 0 \\ 0 & I_{s_v} \end{pmatrix}$.  In particular, if there is $v_0 \in S_\infty$ such that $s_{v_0} > 0$ but $s_v = 0$ for $v \in S_\infty \setminus \{v_0\}$ -- the type of unitary groups that we will be mainly interested in later -- then $E(H,Y_V)$ is the subfield $\imath_{v_0}(F) \subset \CC$. \\\\ 
We will also fix the following notation: Let $V'\subset V$ be a non-degenerate subspace of $V$ of codimension 1. We write $V$ as the orthogonal direct sum $V' \oplus V'_1$ and consider the unitary groups $H':=U(V')$ and $H'':=U(V')\times U(V'_{1})$ over $F^{+}$. Obviously, there are natural inclusions $H'\subset H''\subset H$, and a homomorphism of Shimura data
\begin{equation}\label{inclusionmap} 
(H'',Y_{V'}\times Y_{V'_1}) \hookrightarrow (H,Y_V).
\end{equation}
It is not necessarily the case that $V'_1$, as introduced above, and $V_1$ from \S \ref{sect:alggrp} are isomorphic as hermitian spaces, but the attached unitary groups $U(V'_1)$ and $U(V_1)$ are isomorphic.

 \subsection{Rational structures and (cute) coherent cohomology}\label{sect:RatCute}
\subsubsection{A characterization of cute coherent cohomology}
At each $v\in S_\infty$, we write as usual $\h_{v,\C}= \k_{H,v,\C} \oplus \p^-_v \oplus \p^+_v$ for the Harish-Chandra decomposition of the complex reductive Lie algebra $\h_{v,\C}$, and let
$$\p^+:= \oplus_v  \p^+_v, \quad\quad\p^- := \oplus_v  \p^-_v  \quad\quad{\rm and}\quad\quad \q:=\k_{H,\infty,\C}\oplus \p^-$$
so that
$$\h_{\infty,\C} = \k_{H,\infty,\C} \oplus \p^-  \oplus \p^+=\q\oplus \p^+.$$
Here $\p^+$ and $\p^-$ identify naturally with the holomorphic and anti-holomorphic tangent spaces to $Y_V$ at the point $y$, chosen in order to fix our choice of $K_y=K_{H,\infty}$. The Lie algebra $\q=\q_y$ is a complex parabolic subalgebra of $\h_{\infty,\C}$ with Levi subalgebra $\k_{H,\infty,\C}$. We let $W^\q$ be the set of attached Kostant representatives in the Weyl group of $H_\infty$, cf.\ \cite{bowa}, III.1.4. \\\\
Let  $\lambda = (\lambda_v)_{v \in S_\infty}$ be the highest weight of an irreducible finite-dimensional representation of $H_\infty$ as in \S \ref{sect:coh}. For a $w\in W^\q$, we may form the highest weight $\Lambda(w,\lambda):= w(\lambda+\rho_n)-\rho_n$, $\rho_n$ the half-sum of positive absolute roots of $H_\infty$, of a uniquely determined irreducible, finite-dimensional representation $\W_{\Lambda(w,\lambda)}$ of $K_{H,\infty}$ and we recall that its contragredient $\W^{\sf v}_{\Lambda(w,\lambda)}\cong \W_{\Lambda(w',\lambda^{\sf v})}$ is again of the above form for a uniquely determined Kostant representative $w'\in W^\q$, cf.\ \cite{bowa}, V.1.4. We will henceforth suppress the dependence of $\Lambda$ on $w$ and $\lambda$ in notation.\\\\
Recall from \cite{H90}, \S 2.1, that the representation $\W^{\sf v}_{\Lambda}$ defines an automorphic vector bundle $[\W^{\sf v}_{\Lambda}]$ on the Shimura variety $Sh(H,Y_V)$. Algebraicity of $\lambda$ implies that the canonical and sub-canonical extensions of the $H(\A_{F^+,f})$-homogeneous vector bundle $[\W^{\sf v}_{\Lambda}]$ give rise to coherent cohomology theories which are both defined over a finite extension of the reflex field, see cf.\ \cite{H90}, Proposition 2.8. We let $E(\Lambda)$ denote a number field over which there is such a rational structure.  (In general, there is a Brauer obstruction to realizing $\W^{\sf v}_{\Lambda}$ over the fixed field of its stabilizer in Gal$(\Qbar/\QQ)$, and we can take and fix $E(\Lambda)$ to be some finite, even abelian extension of the latter.)\\\\
Following the notation of \cite{guer-lin}, we denote by $H_!^*([\W^{\sf v}_{\Lambda}])$ the interior cohomology of $[\W^{\sf v}_{\Lambda}]$, cf.\ \cite{H90} \S (3.5.6). This is in contrast to \cite{H14,H90}, where the notation $\bar{H}$ was used.  Interior cohomology, being the image of a rational map, has a natural rational structure over  $E(\Lambda)$. It is well-known that every class in $H_!^*([\W^{\sf v}_{\Lambda}])$ is representable by square-integrable automorphic forms (\cite{H90}, Theorem 5.3) and that the $(\q, K_{H,\infty})$-cohomology of the space of cuspidal automorphic forms injects into $H_!^*([\W^{\sf v}_{\Lambda}])$ (\cite{H90}, Proposition 3.6). Let $H^*_{cute}([\W^{\sf v}_{\Lambda}]) \subseteq H_!^*([\W^{\sf v}_{\Lambda}])$ denote the subspace of classes, represented by cuspidal automorphic forms, contained in {\bf cu}spidal representations that are {\bf te}mpered at all places of $F^+$, where $H$ is unramified. Similarly, let $\CA_{cute}(H)$ be the corresponding space of cuspidal automorphic forms on $H(\A_{F^+})$, which give rise to representations which are tempered at all places of $F^+$, where $H$ is unramified. So, $H^*(\q,K_{ H,\infty},\CA_{cute}(H)\otimes \W^{\sf v}_\Lambda))\cong H^*_{cute}([\W^{\sf v}_{\Lambda}])$.

\begin{prop}\label{h,K-coh}
For a cuspidal automorphic representation $\pi$ of $H(\A_{F^+})$ the following assertions are equivalent:
\begin{enumerate}
\item $\pi \subset \CA_{cute}(H)$ and contributes non-trivially to $H^{*}_{cute}([\W^{\sf v}_{\Lambda}])$ for some $\Lambda=\Lambda(w)$, $w\in W^\q$.
\item $\pi$ is cohomological and its base change $BC(\pi)$ is an isobaric sum $\Pi=\Pi_1\boxplus...\boxplus\Pi_r$ of conjugate self-dual cuspidal automorphic representations $\Pi_i$.
\item $\pi$ is cohomological and tempered.
\end{enumerate}
If $\pi$ satisifies any of the above equivalent conditions, then $\pi_\infty$ is in the discrete series and $\pi$ occurs with multiplicity one in $L^2(H(F^+)\R_+\backslash H(\A_{F^+}))$.
\end{prop}
\begin{proof}
$(1)\Rightarrow (2)$: Let $\pi \subset \CA_{cute}(H)$ denote a cuspidal automorphic representation of $ H(\A_{F^+})$ that contributes to $H^*_{cute}([\W^{\sf v}_{\Lambda}])$. As $\Lambda^{\sf v}=\Lambda(w',\lambda^{\sf v})$ for a (unique) Kostant representative $w'\in W^\q$, it follows from reading the proof of \cite{gro-seb}, Theorem A.1 backwards, that there is an inclusion of vector spaces
$$H^*(\h_\infty, K_{H,\infty}, \pi_\infty \otimes \cF^{\sf v}_\lambda) \hookleftarrow H^*(\q,K_{ H,\infty},\pi_\infty\otimes \W^{\sf v}_\Lambda).$$
Therefore, $\pi_\infty$ is cohomological. Its base change $\Pi=BC(\pi)$ hence exists, cf.\ \S \ref{bc} for our conventions, and is a cohomological isobaric sum $\Pi=\Pi_1\boxplus...\boxplus\Pi_r$ of conjugate self-dual square-integrable automorphic representations $\Pi_i$ of some $G_{n_i}(\A_F)$. As $\pi \subset \CA_{cute}(H)$, $\pi_v$ is unramified and tempered outside a finite set of places $S$ of $F^+$ and hence so is $\Pi$ outside a finite set of places of $F$: Indeed, if $v\notin S$ is split in $F$, then $\Pi_v$ is tempered as noted in \S \ref{bc}. If, however, $v\notin S$ is not split, then $H(F^+_v)\cong U^*_n(F^+_v)$ is the quasisplit unitary group of rank $n$ over $F^+_v$ and so $\pi_v$ has a bounded local Arthur-parameter in the sense of \cite{mok}, Theorem 2.5.1. It follows that the unramfied representation $\Pi_v\cong BC(\pi)_v\cong BC(\pi_v)$ of $G_n(F_v)$ has a bounded local Langlands-parameter, whence $\Pi_v$ is tempered. Now, the argument of the proof of \cite{clozel}, Lemma 1.5, carries over verbatim, showing that the automorphic representation $\Pi$ must be isomorphic to an isobaric sum of unitary cuspidal automorphic representations. By the classification of isobaric sums, cf.\ \cite{jacshal2}, Theorem 4.4, these are nothing else than the isobaric summands $\Pi_i$ from above.

$(2)\Rightarrow (3)$: This is the contents of Remark \ref{rmk:temp}.

$(3)\Rightarrow (1)$: We refer again to \cite{gro-seb}, Theorem A.1, which shows that a cohomological tempered cuspidal automorphic representation has non-trivial $(\q,K_{H,\infty})$-cohomology with respect to a suitable coefficient module $\W^{\sf v}_{\Lambda}$, $\Lambda=\Lambda(w)$, $w\in W^\q$, from which the assertion is obvious. 

In order to prove the last assertions of the proposition, recall from \cite{vozu}, p.\ 58 that a tempered cohomological representation of $H_\infty$, must be in the discrete series. Finally, it follows from Remark 1.7.2 and Theorem 5.0.5 in \cite{KMSW} (and the fact that the continuous $L^2$-spectrum does not contain any automorphic forms) that every $\pi$, which satisfies the equivalent conditions of the proposition, occurs with multiplicity one in $L^2(H(F^+)\R_+\backslash H(\A_{F^+}))$.
 \end{proof}
 
 This result has several consequences. Firstly, we note
 
 \begin{prop}\label{prop:rat}
 The subspace $H^*_{cute}([\W^{\sf v}_{\Lambda}])$ of $H_!^*([\W^{\sf v}_{\Lambda}])$  is rational over $E(\Lambda)$.
\end{prop}
\begin{proof}
Let $H^*_{t}([\W^{\sf v}_{\Lambda}]) \subset H_!^*([\W^{\sf v}_{\Lambda}])$ denote the subspace of interior cohomology, which is represented by forms that are tempered at all non-archimedean places, where the ambient representation is unramified. The condition of temperedness at such a place is equivalent to the condition that the eigenvalues of Frobenius all be $q$-numbers of the same weight, hence is equivariant under Aut$(\C)$. Therefore, $H^*_{t}([\W^{\sf v}_{\Lambda}])$ is an $E(\Lambda)$-rational subspace, and it suffices to show that it coincides with $H^*_{cute}([\W^{\sf v}_{\Lambda}])$. Obviously, by the third item of Proposition \ref{h,K-coh},  $H^*_{cute}([\W^{\sf v}_{\Lambda}])\subseteq H^*_{t}([\W^{\sf v}_{\Lambda}])$, so we may complete the proof by showing that any square-integrable automorphic representation $\pi$ of $H(\A_{F^+})$ that contributes to $H^*_{t}([\W^{\sf v}_{\Lambda}])$ contributes to $H^*_{cute}([\W^{\sf v}_{\Lambda}])$. By \cite{clozel2}, Proposition 4.10, any such $\pi$ must be cuspidal. Now, the argument of the step ``$(1)\Rightarrow (2)$'' of the proof of Proposition \ref{h,K-coh} transfers verbatim, and we obtain that any such $\pi$ satisfies condition (2) of Proposition \ref{h,K-coh}. Hence, again by Proposition \ref{h,K-coh}, $\pi \subset \CA_{cute}(H)$, which shows the claim.
\end{proof}

\begin{rem} 
As far as we know, it has not been proved in general that the cuspidal subspace $H^*(\q,K_{ H,\infty},\CA_{cusp}(H)\otimes \W^{\sf v}_\Lambda))\cong H^*_{cusp}([\W^{\sf v}_{\Lambda}])$ of $H_!^*([\W^{\sf v}_{\Lambda}])$ is rational over $E(\Lambda)$,
but it is known  that for sufficiently regular $\Lambda$ the interior cohomology is entirely cuspidal. In particular, this holds under the regularity assumptions of our main results.
\end{rem}

As another consequence of Proposition \ref{h,K-coh} we obtain

\begin{cor}\label{cor:uniquedegree}
For each $\W^{\sf v}_{\Lambda}$ as in Proposition \ref{h,K-coh}, there is a single degree $q = q(\Lambda)=\sum_{v\in S_\infty}q(\Lambda_v)$, the $q(\Lambda_v)$ being uniquely determined, such that $H^{q(\Lambda)}_{cute}([\W^{\sf v}_{\Lambda}]) \neq 0$. 
\end{cor}
\begin{proof}
Let $v\in S_\infty$. By \cite{H13}, Theorem 2.10, there is a unique discrete series representation $\pi_{\Lambda_v}$ of $H(F_v)$ and a unique degree $q(\Lambda_v)$ such that $H^{q(\Lambda_v)}(\q_v,K_{H,v},\pi_v\otimes \W^{\sf v}_{\Lambda_v})\neq 0$. Moreover, the latter $(\q_v,K_{H,v})$-cohomology is one-dimensional. Hence, by Proposition \ref{h,K-coh}, there are the following isomorphisms for the graded vector space
\begin{equation}\label{dbarcoh}
\begin{aligned}
H^{*}_{cute}([\W^{\sf v}_{\Lambda}]) &\cong \bigoplus_{\pi \subset \CA_{cute}(H)} H^{*}(\q,K_{ H,\infty},\pi_\infty\otimes \W^{\sf v}_\Lambda)\otimes \pi_f \\
&\cong \bigoplus_{\substack{\pi \subset \CA_{cute}( H)\\ \\ \pi_v \simeq \pi_{\Lambda_v}}} \bigotimes_{v\in S_\infty} H^{q(\Lambda_v)}(\q_v,K_{H,v},\pi_v\otimes \W^{\sf v}_{\Lambda_v})\otimes \pi_f \\
&\cong \bigoplus_{\substack{\pi \subset \CA_{cute}( H)\\  \pi_v \simeq \pi_{\Lambda_v} \forall ~v \in S_\infty}} \pi_f.
\end{aligned}
\end{equation}
for the unique degrees $q(\Lambda_v)$. So, $H^{q}_{cute}([\W^{\sf v}_{\Lambda}]) =0$ unless $q= q(\Lambda)=\sum_{v\in S_\infty}q(\Lambda_v)$, in which case $H^{q(\Lambda)}_{cute}([\W^{\sf v}_{\Lambda}]) $ is described by \eqref{dbarcoh}.
\end{proof}

\subsubsection{The field $E(\pi)$}\label{sect:E(pi)}
Let $\pi \subset \CA_{cute}(H)$ be as in the statement of Proposition \ref{h,K-coh}. Recall that the $(\h_\infty, K_{H,\infty},H(\A_{F^+,f}))$-module of smooth and $K_{H,\infty}$-finite vectors in $\pi$ may be defined over a number field $E(\pi)\supseteq\Q(\pi_f)$, cf.\ \cite{H13} Corollary 2.13 \& Proposition 3.17. (Here we use that $\pi$ has multiplicity one in the $L^2$-spectrum, cf.\ Proposition \ref{h,K-coh}, in order to verify the assumption of \cite{H13} Proposition 3.17. See also the erratum to \cite{H13}.) We choose $E(\pi)$ to contain the compositum $F^{Gal}E(\Lambda)$, and refer to these rational structures as the {\it deRham-rational structures} on $\pi$.   These structures are of course inherited from the rational structures on
$H^{q(\Lambda)}_{cute}([\W^{\sf v}_{\Lambda}]) $ for $\Lambda$ and $q(\Lambda)$ uniquely determined by the discrete series $\pi_\infty$.  A function inside this deRham-rational structure is said to be {\it deRham-rational}. We obtain

\begin{lem}\label{contained}
$\Q(BC(\pi)_f^{\sf v})=\Q(BC(\pi)_f)\subseteq E(\pi)$.
\end{lem}
\begin{proof}
Strong multiplicity one implies that $\Q(BC(\pi)_f)=\Q(BC(\pi)^{S})$, where $S$ is any finite set of places containing $S_\infty$ and the places where $BC(\pi)$ ramifies. Hence, $\Q(BC(\pi)_f)=\Q(BC(\pi)^S)=\Q(BC(\pi)^{S,{\sf v}})\subseteq E(\pi^S)$, where the last inclusion is due to \cite{gan-ragh}, Lemma 9.2, the definition of base change and the definition of $E(\pi)$. Invoking strong multiplicity one once more, $\Q(BC(\pi)_f^{\sf v})=\Q(BC(\pi)_f)\subseteq E(\pi)$.
\end{proof}

\begin{lem}\label{lem:unique}
Let $\pi \subset \CA_{cute}(H)$ be an irreducible representation, which contributes non-trivially to $H^{*}_{cute}([\W^{\sf v}_{\Lambda}])$. Then, for each $\sigma\in {\rm Aut}(\C/E(\Lambda))$, there is a unique cohomological tempered cuspidal automorphic representation ${}^\sigma\!\pi$ of $H(\A_{F^+})$, such that $({}^\sigma\!\pi)_f\cong {}^\sigma\!(\pi_f)$ and which contributes non-trivially to $H^{*}_{cute}([\W^{\sf v}_{\Lambda}])$.
\end{lem}
\begin{proof}
Existence follows from Proposition \ref{prop:rat}, Proposition \ref{h,K-coh} and \eqref{dbarcoh}, while uniqueness follows from \cite{H13}, Theorem 2.10, in combination with multiplicity one, see again Proposition \ref{h,K-coh}.
\end{proof}
Our definition of $E(\pi)$ leaves us some freedom to include in it any other appropriate choice of a number field. We will specify such an additional choice right before Conjecture \ref{lvarch}, by adding a suitable number field, constructed and denoted $E_Y(\eta)$ in \cite{H13}, p.\ 2023, to $E(\pi)$. So far, any choice (subject to the above conditions) works.

\subsection{Construction of automorphic motives}\label{construction of motive}

Let $\Pi$ be a cohomological, conjugate self-dual, cuspidal automorphic representation of $\GL_{n}(\Acm)$, which satisfies Hypothesis \ref{descent}. Choose  $I_0=(I_{\imath})_{\imath\in \Sigma} \in \{0,1,\cdots,n\}^{|\Sigma|}$ so that $I_{v_0} = 1$, for some fixed place $v_0$, and so that $I_v = 0$ for $v \neq v_0$. Then, the unitary group $H=H_{I_0}$ has local archimedean signature $(r_{v_0},s_{v_0}) = (n-1,1)$, and the attached group $H'$ from \S \ref{Sdata} above has signature $(r'_{v_0},s'_{v_0}) = (n-2,1)$, while for $v \neq v_0$ in $S_\infty$, the signatures are $(n,0)$ (resp. $(n-1,0)$). \\\\
By Hypothesis \ref{descent}, $\Pi^{\sf v}$ descends to a cohomological tempered cuspidal automorphic representation $\pi$ of $H(\A_{F^+})$. Morover, for any $\pi = \pi_\infty \otimes \pi_f \in \prod(H,\Pi^{\sf v})$ the representation $\tau_\infty\otimes \pi_f$ belongs to $\prod(H,\Pi^{\sf v})$, whenever $\tau_\infty$ is a discrete series representation of $H_\infty$ with the same infinitesimal character as $\pi_\infty$. By Proposition \ref{h,K-coh}, each such $\tau_\infty\otimes \pi_f \in \prod(H,\Pi^{\sf v})$ has multiplicity one in $L^2(H(F^+)\R_+\backslash H(\A_{F^+}))$. (The duality is not a misprint; with the usual normalization it is needed in order to obtain the Galois representation attached to the original $\Pi$, rather than $\Pi^{{\sf v}}$.)\\\\
For a given $\pi_f$ that descends $\Pi_f^\vee$, the set of $\tau_\infty$ such that $\tau_\infty\otimes \pi_f \in \prod(H,\Pi^{\sf v})$ has cardinality $n$, cf.\ \S \ref{sect:coh}. In other words, let $(H,Y_V)$ be the Shimura datum defined in \S \ref{Sdata} and let $Sh(H,Y_V)$ be the corresponding Shimura variety.  There is a unique irreducible finite-dimensional representation $\cF_{\lambda}=\otimes_{v\in S_\infty}\cF_{\lambda_v}$ of $H_\infty$, as in \S \ref{sect:finitereps}, such that, for all $\tau_\infty$ as above, 
\begin{equation}\label{hn-1} 
\dim H^{n-1}(\h_\infty,K_{H,\infty},\tau_\infty\otimes  \cF^{\sf v}_{\lambda}) =1
\end{equation}
Combining this with our previous observations, this implies that 
\begin{equation}\label{motivedimension}  \dim \Hom_{H(\Atrf)}(\pi_f,H^{n-1}(Sh(H,Y_V),\tilde{\cF}^{{\sf v}}_\lambda)) = n
\end{equation}
where $\tilde{\cF}^{{\sf v}}_\lambda$ is the local system on $Sh(H,Y_V)$ attached to the representation $\cF^{{\sf v}}_\lambda$. 
\\\\
The representation $\cF_\lambda$ of $H(F^+)$ is defined over a number field $E(\lambda)$, which we may assume contains the reflex field $E(H,Y_V)=\imath_{v_0}(F)$ of the Shimura variety.  Thus, the cohomology space $H^{n-1}(Sh(H,Y_V),\tilde{\cF}^{{\sf v}}_\lambda)$ has a natural $E(\lambda)$-structure, the {\it Betti} cohomological structure. Letting $\CO(\lambda)$ denote the ring of integers of $E(\lambda)$, we can find a free $\CO(\lambda)$-submodule
$\mathcal M_\lambda \subset \cF_{\lambda}$ that generates the representation, and thus we have a local system in free $\CO(\lambda)$-modules
$$\tilde{\mathcal M}_\lambda^{{\sf v}} \subset \tilde{\cF}^{{\sf v}}_\lambda$$
over $Sh(H,Y_V)$.  
For any prime number $\ell$ and any divisor $\mathfrak{l}$ of $\ell$ in $\CO(\lambda)$ we let 
$$\tilde{\cF}^{{\sf v}}_{\lambda,\mathfrak{l}}:= \tilde{\mathcal M}^{{\sf v}}_\lambda\otimes_{\CO(\lambda)} E(\lambda)_\mathfrak{l} \cong \tilde{\mathcal M}^{{\sf v}}_{\lambda,\mathfrak{l}}\otimes_{\CO(\lambda)_\mathfrak{l}} E(\lambda)_\mathfrak{l}$$
denote the corresponding $\ell$-adic \'etale sheaf.
Then we have the \'etale comparison map
\begin{equation}\label{ladic}
H^{n-1}(Sh(H,Y_V),\tilde{\cF}^{{\sf v}}_{\lambda})\otimes_{E(\lambda)} E(\lambda)_\mathfrak{l} 
 \isoarrow H^{n-1}_{\textrm{{\it \'et}}}(Sh(H,Y_V),\tilde{\cF}^{{\sf v}}_{\lambda,\mathfrak{l}}).
\end{equation}
On the other hand, for any embedding $\imath:  E(\lambda) \hookrightarrow \C$, we have the de Rham comparison 
\begin{equation}\label{dR}
H^{n-1}(Sh(H,Y_V),\tilde{\cF}^{{\sf v}}_{\lambda})\otimes_{E(\lambda),\imath} \C 
 \isoarrow H^{n-1}_{dR}(Sh(H,Y_V),\tilde{\cF}^{{\sf v}}_{\lambda,dR})\otimes_{E(\lambda),\imath} \C.
\end{equation}
Here we let $\tilde{\cF}^{{\sf v}}_{\lambda,dR}$ denote the flat vector bundle over $Sh(H,Y_V)$ attached to the local system $\tilde{\cF}^{{\sf v}}_{\lambda}$ by the Riemann-Hilbert correspondence; the $E(\lambda)$ structure on $H^{n-1}_{dR}$ is derived from the canonical model of $Sh(H,Y_V)$ over $E(H,Y_V) \subset E(\lambda)$, and the rational structure on the flat vector bundle $\tilde{\cF}^{{\sf v}}_{\lambda,dR}$.\\\\  
It is well-known (cf., \cite{H97}, Proposition 2.2.7) that, for any $\lambda$, the Hodge filtration on the right hand side of \eqref{dR} has an associated graded composed of $n$ spaces of interior cohomology:
\begin{equation}\label{Hodge}
gr_F^{\bullet} H^{n-1}_{dR}(Sh(H,Y_V),\tilde{\cF}^{{\sf v}}_{\lambda,dR}) \cong\bigoplus_{q = 0}^{n-1} H_!^q([\W^{\sf v}_{\Lambda(q)}]),
\end{equation}
Here $\Lambda(q) = \HC(q) - \rho_n$, where $\HC(q)$ is the Harish-Chandra parameter of $\pi_{\lambda,q}$, cf.\ \S \ref{sect:coh}. For $q = 0, \dots, n-1$, we define $i(q) \in \Z$ by
$$H_!^q([\W^{\sf v}_{\Lambda(q)}])= gr_F^{i(q)}H^{n-1}_{dR}(Sh(H,Y_V),\tilde{\cF}^{{\sf v}}_{\lambda,dR}).$$
We now recall that the group $H(\Atrf)$ acts on the spaces in \eqref{ladic}, \eqref{dR}, and \eqref{Hodge} compatibly with the comparison isomorphisms.  Let $\pi_f$ be as before. It is defined over the number field $E(\pi) \supset E(\lambda)$, as introduced in \S \ref{sect:E(pi)}, and we define
$$M_{dR}(\pi_f) := \Hom_{H(\Atrf)}(\pi_f,H_{dR}^{n-1}(Sh(H,Y_V),\tilde{\cF}^{{\sf v}}_{\lambda,dR})\otimes E(\pi)).$$
$$M_B(\pi_f):= \Hom_{H(\Atrf)}(\pi_f,H^{n-1}(Sh(H,Y_V),\tilde{\cF}^{{\sf v}}_{\lambda})\otimes E(\pi)),$$
$$M_{\mathfrak{l}}(\pi_f) := \Hom_{H(\Atrf)}(\pi_f,H^{n-1}_{\textrm{{\it \'et}}}(Sh(H,Y_V),\tilde{\cF}^{{\sf v}}_{\lambda,\mathfrak{l}}).),$$
Here we are abusing notation: The $\pi_f$ in each $\Hom$ space above is viewed as a vector space over the appropriate coefficient field by extension of scalars, namely $E(\pi)$, in the first two, and $E(\pi)_{\mathfrak{l}}$, in the third line. Clearly all three of the spaces $M_?(\pi_f)$ have the same dimension over their respective coefficient fields.  In fact, one may check that $\dim_{E(\pi)} M_B(\pi_f) = n$.  More precisely, the Hodge filtration on $H^{n-1}_{dR}(Sh(H,Y_V),\tilde{\cF}^{{\sf v}}_{\lambda,dR})$ induces a decreasing filtration $F^{i}M_{dR}(\pi_f)$ on $M_{dR}(\pi_f)$, and the isomorphism $\eqref{Hodge}$ induces an isomorphism
$$gr_F^{\bullet}M_{dR}(\pi_f) = \bigoplus_{q=0}^{n-1}gr_F^{i(q)} M_{dR}(\pi_f)  \cong \bigoplus_{q = 0}^{n-1} \Hom_{H(\Atrf)}(\pi_f,H_!^q([\W^{\sf v}_{\Lambda(q)}])), $$
and each of the spaces $M^{i(q)}_{dR}(\pi_f):= gr_F^{i(q)} M_{dR}(\pi_f)$ is of dimension $1$ over $E(\pi)$. The following result now follows from our construction:

\begin{thm}\label{thm:clozel}
The collection $(M_B(\pi_f), M_{dR}(\pi_f), \{M_{\mathfrak{l}}(\pi_f)\}_{\mathfrak l})$, together with the obvious comparison maps defines a regular, pure motive $M(\Pi)$ over $E(H,Y_V)$ with coefficients in the finite extension $E(\pi)$ of $\Q(\Pi_f)$.  More precisely, the data satisfy the conditions of Definition \ref{definitionmotive}, with the exception of (i) (the infinite Frobenius); see \ref{Finfty} below.

Moreover, if $\Pi'$ is another cohomological, conjugate self-dual, cuspidal automorphic representation of $\GL_{n'}(\Acm)$, which satisfies Hypothesis \ref{descent}, then 
$$L(s,M(\Pi)\otimes M(\Pi'))=L(s-\tfrac{n+n'-2}{2},\Pi_f\times \Pi'_f),$$
interpreted as $E(\pi)E(\pi')\otimes_\Q\C$-valued functions as in \S\ref{tensorp}.
\end{thm}
 
A few remarks are in order:

\begin{rmk}  
(1) If we recall that $E(H,Y_V)\cong F$, putting $n'=1$ in Theorem \ref{thm:clozel} proves \cite{clozel}, Conjecture 4.5, for the conjugate self-dual, cuspidal automorphic representations at hand.\\\\
(2) The motive $M(\Pi)$ depends on the choice of the place $v_0$, as so does the Shimura variety $Sh(H,Y_V)$. In view of \S\ref{sect:alggrp} and \cite{milne-suh}, Theorem 1.3, replacing $v_0$ by a different choice $v_1$ means to descend to the unitary group ${}^\sigma H$ underlying the $\sigma$-twisted Shimura variety ${}^\sigma Sh(H,Y_V)$, where $\sigma$ is any complex automorphism such that $\sigma^{-1}\circ \imath_{v_0} = \imath_{v_1}$. Hence, upon applying restriction of scalars, one obtains a motive $R_{E(H,Y_V)/\Q}(M(\Pi))$ over $\Q$, which is in fact independent of the choice of $v_0$.\\\\
(3) In the setting of the PEL type Shimura variety attached to the unitary similitude group containing the unitary group $H$, Theorem \ref{thm:clozel} follows from the construction of Galois representations attached to cohomological cuspidal representations of $\GL_n$ over the CM field $F$, starting with \cite{clozelihes} and continuing over more than 20 years through \cite{car}.  In particular, the identification of the $L$-functions of tensor products of motives with the Rankin-Selberg automorphic $L$-functions makes use of the local Langlands correspondence in the form proved in \cite{HT}.  

The article \cite{KSZ} carries out the analysis of points on special fibers of Shimura varieties of abelian type, including $Sh(H,Y_V)$. The applications to the global $L$-functions has not yet been written up.  The cautious reader may prefer to consider the statements in this section to be special cases of the Langlands conjecture on Hasse-Weil zeta functions of Shimura varieties, whose proof, at least at unramified places, has been promised for a sequel to \cite{KSZ}. Since the statements claimed here will not be used elsewhere, except heuristically in the statement of Theorem \ref{main factorization}, this is harmless.  Moreover, the $L$-functions in the statement of Theorem \ref{thm:clozel} can be replaced by the partial $L$-functions, with ramified factors removed, with no effect on the rationality results of the present paper.
\end{rmk} 

At least if $F^+ \neq \Q$, the motive $M(\Pi)$ can be identified with a direct summand in the cohomology of a certain abelian scheme over a locally symmetric space $S'(H)$ isomorphic over $\Qbar$ to $Sh(H,Y_V)$ -- but with algebraic structure inherited from an embedding in the PEL Shimura variety attached to a similitude group containing $H$. 

\subsubsection{The automorphic version of $F_\infty$}\label{Finfty}  It is most convenient to take complex conjugation of differential forms as a surrogate for the operator $F_{B,\imath}$ of Definition \ref{definitionmotive}.  For the reason explained in \cite[Remark 3.5]{H21}, this is not quite right.  This is why the automorphic $Q$-periods of \S \ref{Qperiod}, which arise naturally in the calculation of $L$-functions, do not quite correspond to the motivic periods of \S \ref{motivic periods}.  We return to this point in \S \ref{Qperiod} and in \S \ref{local period definition}.

\section{Periods for unitary groups and the Ichino-Ikeda-Neal Harris conjecture}

\subsection{GGP-periods, pairings for unitary groups, and a recent theorem}\label{sect:GGPunit}

Let $V$, $V'$, $V'_1$, $H$, $H'$, $H''$ be as in \S \ref{Sdata}. The usual Ichino-Ikeda-N. Harris conjecture considers the inclusion $H'\subset H$. However, in view of \eqref{inclusionmap} it is sometimes more convenient to consider the inclusion $H''\subset H$ instead, see \cite{H13} and \cite{H14}, and we are going to use both points of view in this paper. In this section we take the opportunity to discuss the relations of the associated periods for the two inclusions $H'\subset H$ and $H''\subset H$. We warn the reader that our notation here differs slightly from \cite{H13} and \cite{H14}.\\\\
Let $\pi$ (resp. $\pi'$) be a cohomological tempered cuspidal automorphic representation of $H(\Atr)$ (resp.\ $H'(\Atr)$). Let $\xi$ be a Hecke character on $U(V_{1})(\Atr)$ (recall that $U(V_1)$ is independent of the hermitian structure on $V_1$, \S \ref{sect:alggrp}). We write $\pi'':=\pi' \otimes \xi$, which is a tempered cuspidal automorphic representation of $H''(\Atr)$. Moreover, we fix a Haar measure $dh:=\prod_{v} dh_v$ on $H(\Atr)$, normalized as in \cite{H13}, \S 5.2, adding the (compatible) convention that $vol_{dh}(U(V_1)(F^+)\backslash U(V_1)(\A_{F^+}))=1$. This defines measures on $H'(\Atr)$ and $H''(\Atr)$ accordingly.\\\\
For  $f_{1}, f_{2}\in \pi$ the {\it Petersson inner product} on $\pi$ is defined as usual as
$$ \<f_{1},f_{2}\>:=\int_{H(\tr)Z_{H}(\Atr)\backslash H(\Atr)} f_{1}(h)\overline{f_{2}(h)} \ dh.$$
Analogously, we may define the Petersson inner product on $\pi'$ and $\pi''$. Next, for $f\in \pi$, $f'\in \pi'$ we put
$$I^{can}(f,f'):=\int\limits_{H'(\tr)\backslash H'(\Atr)} f(h')f'(h') \ dh' .$$
Then it is easy to see that $I^{can} \in \Hom_{H'(\Atr)}(\pi\otimes \pi',\C)$. With this notation the {\it GGP-period} for the pair $(\pi,\pi')$ 
(called the Gross-Prasad period in \cite{H13}) is defined as 
$$\CP(f,f'):=\cfrac{|I^{can}(f,f')|^{2}}{\<f,f\>\ \<f',f'\>}.$$
We can similarly define a $H''(\Atr)$-invariant linear form on $\pi \otimes \pi''$, which we will also denote by $I^{can}$, as
$$I^{can}(f,f''):=\int\limits_{H''(\tr)\backslash H''(\Atr)} f(h'')f''(h'') \ dh'' \text{ for } f\in \pi, f''\in \pi'',$$
leading to a definition of the {\it GGP-period} for the pair $(\pi,\pi'')$ as 
$$\CP(f,f''):=\cfrac{|I^{can}(f,f'')|^{2}}{\<f,f\>\ \<f'',f''\>}.$$
Let $\xi_\pi$, resp. let $\xi_{\pi'}$, be the central character of $\pi$, resp.\ $\pi'$. We assume that 

\begin{equation}\label{central extra U(1)}
\xi_{\pi}^{-1}=\xi_{\pi'}\xi 
\end{equation} 
(resembling equation $(\Xi)$ on page 2039 of \cite{H13}). Then, without restriction of generality, we may write a $f'' \in \pi''$, as $f''=f'\cdot \xi$, and one verifies easily that with our normalizations
$$I^{can}(f,f'')=I^{can}(f,f''|_{H'(\Atr)})=I^{can}(f,f')$$ 
and 
$$\<f'',f''\>=\<f''|_{H'(\Atr)},f''|_{H'(\Atr)}\>= \<f',f'\>.$$ We conclude that:
\begin{lem}\label{add U1}
$$\CP(f,f'')=\CP(f,f''|_{H'(\Atr)})=\CP(f,f').$$
Moreover, one gets $\Hom_{H'(\Atr)}(\pi\otimes \pi',\C) \cong \Hom_{H''(\Atr)}(\pi\otimes \pi'',\C)$.
\end{lem}
We will also need a local version of the above pairings. To this end, choose $f \in \pi$, $f' \in \pi'$, and assume they are factorizable as
$f = \otimes f_v, f' = \otimes f'_v$ with respect to the restricted tensor product factorizations
\begin{equation}\label{422} 
\pi \cong \otimes'_v \pi_v, \quad\quad \pi' \cong \otimes'_v \pi'_v.  
\end{equation}
Outside a finite set $S\supset S_\infty$ of places of $F^+$, we assume $\pi_v$ and $\pi'_v$ are unramified, and $f_v$ and $f'_v$ are the normalized
spherical vectors, i.e., the unique spherical vector taking value $1$ at the identity element.  We choose inner products $\<\cdot,\cdot\>_{\pi_v}, \<\cdot,\cdot\>_{\pi'_v}$
on each of the unitary representations $\pi_v$ and $\pi'_v$ such that at an unramified place $v$, the local normalized spherical vector in $\pi_v$ or $\pi'_v$ has norm $1$. For each place $v$ of $F^+$, let
$$c_{f_v}(h_v) := \<\pi_v(h_v) f_v,f_v\>_{\pi_v} \quad\quad\quad c_{f'_v}(h'_v) := \<\pi'_v(h'_v) f'_v,f'_v\>_{\pi'_v},  ~h_v \in H_v, h'_v \in H'_v,$$
and define
$$I_v(f_v,f'_v) := \int_{H'_v} c_{f_v}(h'_v) c_{f'_v}(h'_v) dh'_v  \quad\quad\quad I^*_v(f_v,f'_v) := \frac{I_v(f_v,f'_v)}{c_{f_v}(1)c_{f'_v}(1)}.$$
Neal Harris proves that these integrals converge since $\pi$ and $\pi'$ are locally tempered at all places.\\\\ 
The GGP-periods and local pairings are interconnected by the Ichino-Ikeda-N.Harris conjecture, which is now a theorem: In order to state it, denote the base change of the cohomological  tempered cuspidal automorphic representation  $\pi\otimes\pi'$ of  $H(\A_{F^+})\times H'(\A_{F^+})$ to $G_n(\A_{F})\times G_{n-1}(\A_{F})$ by $\Pi\otimes\Pi'$. We define
\begin{equation}\label{421} \CL^S(\Pi,\Pi'): = \frac{L^S(\tfrac{1}{2},\Pi\otimes \Pi')}{L^S(1,\Pi,{\rm As}^{(-1)^n})L^S(1,\Pi',{\rm As}^{(-1)^{n-1}})},\end{equation}
where $L^S(1,\Pi,{\rm As}^\pm)$ denotes the partial Asai $L$-function of the appropriate sign and we let
$$\Delta_{H}:= \prod_{i = 1}^n L(i, \varepsilon^i)$$ 
The Ichino-Ikeda-N.Harris conjecture for unitary groups is now the following theorem

\begin{thm} \label{conjecture II} 
Let $f \in \pi$, $f' \in \pi'$ be factorizable vectors as above. Then there is an integer $\beta$ (depending on the Arthur-Vogan-packets containing $\pi$ and $\pi'$), such that
$$\CP(f,f') = 2^{\beta}\Delta^S_{H} \CL^S(\pi,\pi') \ \prod_{v \in S} I^*_v(f_v,f'_v) .$$
\end{thm}

\begin{rem}\label{rem:IIproved}
The conjecture has been proved in increasingly general versions in \cite{zhang, xue, BP3}, and finally the proof was completed in \cite{BLZZ19,BCZ20}.  For totally definite unitary groups it was also shown in \cite{grob_lin} up to a certain algebraic number, under the assumption that the base change of $\pi$ to $G_n(\A_F)$ is cuspidal. 
\end{rem}

\begin{rem}\label{Ivstar}  
Both sides of Theorem \ref{conjecture II} depend on the choice of factorizable vectors $f, f'$, but the dependence is invariant under scaling.  In particular, the statement is independent of the choice of factorizations \eqref{422}, and the assertions below on the nature of the local factors $I^*_v(f_v,f'_v)$ are meaningful.
\end{rem}

The algebraicity of local terms $I^*_v$ was proved in \cite{harrisANT} when $v$ is non-archimedean.  More precisely, we have the following:

\begin{lem}\label{localIv} 
Let $v$ be a non-archimedean place of $F^+$. Let $\pi$ and $\pi'$ be cohomological tempered cuspidal automorphic representations as above. Let $E$ be a number field over which $\pi_v$ and $\pi'_v$  both have rational models.  Then for any $E$-rational vectors $f_v \in \pi_v$, $f'_v \in \pi'_v$, we have
$$I^*_v(f_v,f'_v)  \in E.$$
Moreover, for all $\sigma \in Aut(\CC)$,
$$\sigma\left(I_{v}^{*}\left(f_{v}, f_{v}^{\prime}\right)\right)=I_{v}^{*}\left({ }^{\sigma} f_{v},{ }^{\sigma} f_{v}^{\prime}\right)$$
\end{lem}
\begin{proof}  The algebraicity of the local zeta integrals $I_v(f_v,f'_v)$ is proved in \cite{harrisANT}, Lemma 4.1.9, when the local inner products $\<\cdot,\cdot\>_{\pi_v}$ and  $\<\cdot,\cdot\>_{\pi'_v}$ are taken to be rational over $E$.  Since $f_v$ and $f'_v$ are $E$-rational vectors, this implies the assertion for the normalized integrals $I^*_v$ as well.   The proof
of the first assertion in \cite{harrisANT} is based on the analysis of the local integrals by Moeglin and Waldspurger in \cite{MW12}.  The same analysis shows that the normalized integrals are equivariant with respect to the action of $Aut(\CC)$, and thus implies the second assertion.
\end{proof}

We will state an analogous result for the archimedean local factors as an expectation of ours in the next section.

\subsection{Review of the results of \cite{H14}}\label{reviewH14}

Let us fix a place $v_0\in S_\infty$ and let $H=H_{I_0}$, where $I_0$ is as in \S \ref{construction of motive}. Let $\cF_{\lambda}=\otimes_{v\in S_\infty}\cF_{\lambda_v}$ be some irreducible finite-dimensional representation of $H_\infty=H_{I_0,\infty}$ as in \S \ref{sect:finitereps}. Using that $H_\infty$ is compact at $v\neq v_0$ one sees as in \S \ref{sect:coh} that there are exactly $n$ inequivalent discrete series representations of $ H_{\infty}$, denoted $\pi_{\lambda,q}$, $0\leq q \leq n-1$, for which  $H^p(\h_\infty, K_{H,\infty},\pi_{\lambda,q}\otimes\cF^{\sf v}_{\lambda}) \neq 0$ for some degree $p$ (which necessarily equals $p=n-1$). Moreover, the representations $\pi_{\lambda,q}$, $0\leq q \leq n-1$, are distinguished by the property that,
$$\dim H^q(\q,K_{H,\infty},\pi_{\lambda,q}\otimes \W^{\sf v}_{\Lambda(q)}) = 1$$
for $\Lambda(q)=A(q)-\rho_n$, where $A(q)$ denotes the Harish-Chandra parameter of $\pi_{\lambda,q}$ and all other $H^*(\q,K_{H,\infty},\pi_{\lambda,q}\otimes \W)$ vanish as $\W$ runs over all irreducible representations of $K_{H,\infty}$. We can determine $\HC(q)$ explicitly: Let
$$\HC_{\lambda} = (\HC_{\lambda,v})_{ v \in S_\infty};  ~\HC_{\lambda,v} =(\HC_{v,1} > \dots > \HC_{v,n})$$
be the infinitesimal character of $\cF_{\lambda}$, as in \S 4.2 of \cite{H14}.  Then $\HC(q) = (\HC(q)_{v})_{ v \in S_\infty}$ where $\HC(q)_{v} = \HC_{\lambda,v}$ for $v \neq v_0$ and
$$\HC(q)_{v_0} = (\HC_{v_{0},1} > \dots > \widehat{\HC_{v_{0},q+1}}> \dots  > \HC_{v_{0},n}; \HC_{v_{0},q+1})$$
(the parameter marked by $\widehat{}$ is deleted from the list).  The following is obvious:

\begin{lem}\label{knownknown}  
For $0\leq q \leq n-2$ the parameter $\HC(q)$ satisfies Hypothesis 4.8 of \cite{H14}. For $q=0$, the representation $\pi_{\lambda,q}$ is holomorphic.
\end{lem}

Now suppose that the highest weight $\lambda_{v_0}$ is regular. Equivalently, the Harish-Chandra parameter $\HC_{\lambda,v_0}$ satisfies the regularity condition $\HC_{v_0,i} - \HC_{v_0,i+1} \geq 2$ for $i = 1, \dots, n-1$. Then, for $0\leq q \leq n-2$ define a Harish-Chandra parameter $\HC'(q) = (\HC'(q)_v)_{ v \in S_\infty}$ by the formula (4.5) of \cite{H14}:
\begin{equation}\label{lambdaprime}  
\HC'(q)_{v_0} = (\HC_{v_{0},1} - \tfrac{1}{2} > \dots > \widehat{\HC_{v_{0},q+1}-\tfrac{1}{2}} > \dots  > \HC_{v_{0},n-1}  - \tfrac{1}{2}; \HC_{v_{0},q+1} +  \tfrac{1}{2}).
\end{equation}
For $v \neq v_0$, $\HC'(q)_v =  (\HC_{v,1} - \tfrac{1}{2}  > \dots  > \HC_{v,n-1}  - \tfrac{1}{2})$.\\\\
Since $\lambda_{v_0}$ is regular, \cite{H14}, Lemma 4.7, shows that $\HC'(q)$ is the Harish-Chandra parameter for a unique discrete series representation $\pi_{\HC'(q)}$ of $H'_{\infty}$;
we define the discrete series representation $\pi_{A(q)}$ of $H_\infty$ analogously.  Indeed, the regularity of $\lambda_{v_0}$ is the version of Hypothesis 4.6 of \cite{H14}, where the condition is imposed only at the place $v_0$ where the local unitary group is indefinite. Observe that there is no need for a regularity condition at the definite places: For $v \neq v_0$ the parameter $A'(q)_v$ is automatically the Harish-Chandra parameter of an irreducible representation. We can thus adapt Theorem 4.12 of \cite{H14} to the notation of the present paper:

\begin{thm}\label{GPknownknown}  
Suppose $\lambda_{v_0}$ is regular.  For $0\leq q \leq n-2$ let $ \pi(q)= \pi(A(q))$ and $ \pi'(q)= \pi'(A'(q))$ be a tempered cuspidal automorphic representations of $ H(\A_{F^+})$ and $ H'(\A_{F^+})$, respectively, with archimedean components $\pi_{\HC(q)}$ and $\pi^{\sf v}_{\HC'(q)}$. Let $ \xi$ be the Hecke character  $(\xi_{ \pi}\cdot\xi_{ \pi'})^{-1} $ of $U(V_{1})(\A_{F^+})$ and set $\pi''(q)=\pi'(q)\otimes\xi$. Then for any deRham-rational elements $f\in \pi(q)$, $f''\in \pi''(q)$
$$I^{can}(f,f'')\in E(\pi(q))  E(\pi''(q))=E(\pi(q))  E(\pi'(q)).$$ 
\end{thm}
\noindent The statement in \cite{H14} has two hypotheses:  the first one is the regularity of the highest weight, while the second one (Hypothesis $4.8$ of \cite{H14}) follows as in Lemma \ref{knownknown} from the assumption that $q \neq n-1$. We remark that the assumption in  \textit{loc.cit} on the Gan-Gross-Prasad multiplicity one conjecture for real unitary groups has been proved by He in \cite{he}.\\\\
A cuspidal automorphic representation $ \pi$ of $ H(\A_{F^+})$ that satisfies the hypotheses of Theorem \ref{GPknownknown} contributes to interior cohomology of the corresponding Shimura variety $Sh(H, Y_V)$ with coefficients in the local system defined by the representation $\W^{\sf v}_{\Lambda(q)}$. This cohomology carries a (pure) Hodge structure of weight $n-1$, with Hodge types corresponding to the infinitesimal character of $\W^{\sf v}_{\Lambda(q)}$, which is given by 
$$A(q)^{{\sf v}}= (-A(q)_{v,n}>...> -A(q)_{v,1})_{v\in S_\infty}.$$
The Hodge numbers corresponding to the place $v_0$ are
\begin{equation}\label{hodge1} 
(p_i = -\HC_{v_0,n+1 -i}+ \tfrac{n-1}{2}, q_i = n-1-p_i);  \quad\quad (p_i^c = q_i, q_i^c = p_i). 
\end{equation}
Analogously, a cuspidal automorphic representation $\pi'$ of $ H'(\A_{F^+})$ as in Theorem \ref{GPknownknown} contributes to interior cohomology of the corresponding Shimura variety $Sh(H',Y_{V'})$  with coefficients in the local system defined by a representation $\W_{\Lambda'(q)}$, whose parameters are obtained from those of $\HC'(q)$ by placing them in decreasing order and substracting $\rho_{n-1}$.  In particular, it follows from \eqref{lambdaprime} that the infinitesimal character of $\W_{\Lambda'(q)}$ at $v_0$ is given by
$$(\HC_{v_0,1} - \tfrac{1}{2} > \dots > \HC_{v_0,q} - \tfrac{1}{2} > \HC_{v_0,q+1} + \tfrac{1}{2} > \HC_{v_0,q+2} - \tfrac{1}{2} > \dots > \HC_{v_0,n-1}  - \tfrac{1}{2}),$$
with strict inequalities due to the regularity of $A_{\lambda,v_0}$ and with corresponding Hodge numbers
\begin{equation}\label{hodge2} p_i' = \HC_{v_0,i} + \tfrac{n-3}{2}, \quad {\rm for}\quad i \neq q+1; \quad\quad p'_{q+1} = \HC_{v_0,q+1} + \tfrac{n-1}{2}
 \end{equation}
and $q'_i = n-2 - p'_i$, etc. Here is a consequence of the main result of \cite{H14}.

\begin{thm}\label{BS}  
Let ${\pi}(q)$ be as above. Then -- up to possibly replacing $H'$ by an inner form with the same signatures at all archimedean places -- there exists a tempered cuspidal automorphic representation $\pi'(q)$ of $H'(\A_{F^+})$ with archimedean component $\pi_{\HC'(q)}^{{\sf v}}$, such that 
\begin{enumerate}
\item\label{BS1} $BC(\pi'(q))$ is cuspidal automorphic and supercuspidal at a non-archimedean place of $\tr$ which is not split in $F$, and
\item\label{BS2} there are factorizable cuspidal automorphic forms $f \in \pi(q)$, $f' \in \pi'(q)$, so that $I^{can}(f,f') \neq 0$ with $f_v$ (resp.\ $f'_v$) in the minimal $K_{H,v}$- (resp.\ $K_{H',v}$-type) of $\pi(q)_v$  (resp.\ $\pi'(q)_v$) for all $v \in S_{\infty}$. 
\end{enumerate}
In particular, the GGP-period $\CP(f,f')$ does not vanish.
\end{thm}

\begin{proof}  
Although this is effectively the main result of \cite{H14}, it is unfortunately nowhere stated in that paper. So, let us explain why this is a consequence of the results proved there. First, we claim that the discrete series representation $\pi_{\HC'(q)}$ is isolated in the (classical) automorphic spectrum of $H'_\infty$, in the sense of \cite{BS}, see Corollary 1.3 of \cite{H14}. Admitting the claim, we note that Hypothesis 4.6 of \cite{H14} is our regularity hypothesis, and Hypothesis 4.8 is true by construction. The theorem then follows from the discussion following the proof of Theorem 4.12 of \cite{H14}. More precisely,  because $\pi_{\HC'(q)}$ is isolated in the automorphic spectrum, we can apply Corollary 1.3 (b) of \cite{H14}.  As explained in the discussion of \cite{HLi}, Proposition 3.1, this is a restatement of the main result of \cite{BS}.\\\\
For condition (1), we actually need the $S$-arithmetic version of the Burger-Sarnak result, proved in \cite{CU04}, Theorem 5.3.  We need to show that for some non-archimedean place $w$ of $\tr$, not split in $F$, there is a representation $\tau$  of $H'(\tr_w)$ such that $BC(\tau)$ is supercuspidal and such that 
\begin{equation}\label{ggploc} {\rm Hom}_{H'(\tr_w)}(\pi(q)_w\otimes \tau,\CC) \neq 0. \end{equation}
Here $\tau$ is any member of the $L$-packet whose base change is $BC(\tau)$; but since the base change is supercuspidal it is known that the $L$-packet is
a singleton.\\\\
Let 
$$\sigma:  \Gamma_{F_w} := {\rm Gal}(\bar{F}_w/F_w)  \ra {\rm GL}_n(\CC)$$
denote the Galois parameter that corresponds to $BC(\tau)$ under the local Langlands correspondence.  The requirement is that $\sigma$ be irreducible.
By the local Langlands correspondence for unitary groups, established in \cite{KMSW}, $\tau$ must then also have an irreducible Langlands parameter, and thus is also supercuspidal.   For $w$ of sufficiently large residue characteristic  $H(\tr_w)$ and $H'(\tr_w)$ are quasi-split and $\pi(q)_w$ is unramified; we assume this to be the case.  
It thus suffices to find a supercuspidal $\tau$ of (the quasi-split) $H'(\tr_w)$ satisfying
\eqref{ggploc} for an unramified $\pi(q)_w$, such that $BC(\tau)$ is supercuspidal.    If $n$ is even then $\pi(q)_w$ does not transfer to the non-quasi-split form of $H$ over $w$, and then any supercuspidal $\tau$ satisfies \eqref{ggploc} by the Gan-Gross-Prasad Conjecture for unitary groups \cite{BP2}.  
If $n$ is odd then any supercuspidal representation $\tilde{\tau}$ of ${\rm GL}_{n-1}(F_w)$ descends both to $H'(\tr_w)$ and to its non-quasi-split inner form 
$H''$ over $\tr_w$, necessarily to a supercuspidal representation.  Again by the Gan-Gross-Prasad Conjecture, \eqref{ggploc} holds, up to replacing $H'(\tr_w)$
by $H''(\tr_w)$.  For any $n$, therefore, we may thus choose a global inner form of $H'$ satisfying \eqref{ggploc} at $w$; this may require changing the inner form at some other non-archimedean place but we do not change the signatures at archimedean primes.\footnote{We thank Dipendra Prasad for convincing us that, even when both $H$ and $H'$ are quasi-split at $w$, there is no obvious way to complete this step without invoking the Gan-Gross-Prasad Conjecture.}\\\\
It thus suffices to show that there is some supercuspidal representation of ${\rm GL}_{n-1}(F_w)$ that descends to the group $H'(\tr_w)$; in other words, that there is a stable supercuspidal $L$-packet for $H'(\tr_w)$.  Equivalently, by the local Langlands correspondence, it suffices to show that there is an irreducible parameter $\sigma$ 
whose ${\rm Gal}(F_w/\tr_w)$-conjugate is isomorphic to $\sigma^\vee$.   We may choose $w$ so that $\tr_w$ has a totally ramified cyclic extension $K$ of degree $n-1$.
Let $L = K\cdot F_w$, an abelian extension of $\tr_w$ with Galois group $\Gamma = {\rm Gal}(F_w/\tr_w)\times {\rm Gal}(K/\tr_w)$.   Let $c \in {\rm Gal}(F_w/\tr_w)$ be the non-trivial element.\\\\
Let $p$ be the residue characteristic of $F_w$ and let $W^+ \subset \CO_L^\times$ be the group of units congruent to $1$ modulo the maximal ideal of $\CO_L$.  Then $W = W^+ \otimes_{\Z_p} \QQ_p$ is a free 
$\QQ_p[\Gamma]$-module.  There are thus homomorphisms $\lambda:  W \ra \QQ_p$ such that 
 $\lambda^c = -\lambda$ and 
 $\lambda$ is stabilized by no non-trivial element of ${\rm Gal}(K/\tr_w) = {\rm Gal}(L/F_w)$. By composing such a $\lambda$ with an appropriate character $\psi:  \QQ_p \ra \CC^\times$, we can construct a character $\chi:  W^+ \ra \CC^\times$
such that
\begin{enumerate}
\item $\chi^c = \chi^{-1}$ ;
\item $\chi$ is stabilized by no non-trivial element of ${\rm Gal}(L/F_w)$.
\end{enumerate}  
We can view $\chi$ as a totally ramified
character of $\Gamma_L: = {\rm Gal}(\bar{L}/L)$.  Taking $\sigma$ to be the induced representation $I_{\Gamma_L}^{\Gamma_{F_w}}(\chi)$, (2) implies that $\sigma$ is irreducible and (1) implies that $\sigma^c \cong \sigma^\vee$. Finally, the isolation follows as in \cite{HLi}, Theorem 7.2.1, using the existence of base change from $H'(\A_{F^+})$ to $G_{n-1}(\A_F)$, as was established for tempered cuspidal automorphic representations in \cite{KMSW}, Theorem 5.0.5.
\end{proof}

Recall the number field $E_Y(\eta)$ from \cite{H13}, p.\ 2023. It has been shown in Corollary 3.8 of \cite{H13} (and its correction in the Erratum to that paper) that the underlying Harish-Chandra modules of the discrete series representations $\pi(q)_v$ and $\pi(q)'_v$, $v\in S_\infty$, are defined over this number field $E_Y(\eta)$. {\it From now on, we will assume that the number field $E(\pi)$, defined for a cohomological tempered cuspidal representation of a unitary group over $F^+$ in \S \ref{sect:E(pi)}, contains $E_Y(\eta)$.} One sees that the cuspidal automorphic forms $f \in \pi(q)$, $f' \in \pi'(q)$ from Theorem \ref{BS}.(2) can be chosen so that, for all $v\in S_\infty$, $f_v$ (resp. $f'_v$) belongs to the $E(\pi(q))$- (resp. $E(\pi'(q))$-) rational subspace of the minimal $K_{H,v}$-type of $\pi(q)_v$ (resp. $K_{H',v}$-type of $\pi'(q)_v$), with respect to the $E(\pi(q))$- (resp. $E(\pi'(q))$-) deRham-rational structure defined in \S \ref{sect:E(pi)}. The following statement is then a conjectural, archimedean analog of Lemma \ref{localIv}:

\begin{conj}\label{lvarch}  
Let $f$ and $f'$ be as in Theorem \ref{BS} and assume that they are chosen so that, for all $v\in S_\infty$, $f_v$ (resp. $f'_v$) belongs to the $E(\pi(q))$- (resp. $E(\pi'(q))$-) rational subspace of the minimal $K_{H,v}$-type of $\pi(q)_v$ (resp. $K_{H',v}$-type of $\pi'(q)_v$). Then,
$$I^*_v(f_v,f'_v) \in (E(\pi(q))\cdot E(\pi'(q)))$$
for all $v\in S_\infty$.

\end{conj}

\subsection{Automorphic $Q$-periods}\label{Qperiod}
In order to prove the factorization of the local arithmetic automorphic periods, see \eqref{eq:factor}, we will need one last ingredient, namely {\it automorphic $Q$-periods}. To define them, let $\pi$ be a cohomological, tempered, cuspidal automorphic representation of $H(\A_{F^+})$ and let $\varphi_\pi\in \pi$ be deRham-rational, cf.\ \S \ref{sect:E(pi)}. For each $\sigma\in {\rm Aut}(\C/E(\Lambda))$, we choose a ${}^{\sigma}\! \varphi_\pi$, which generates the deRham-rational structure of the unique twist ${}^\sigma\!\pi$, cf.\ Lemma \ref{lem:unique}. We define
$$Q({}^{\sigma}\!\varphi_{\pi}):=\<{}^{\sigma}\!\varphi_{\pi},{}^{\sigma}\!\varphi_{\pi}\>.$$ 
By Lemma 3.19 and Lemma 3.20 of \cite{H13}, $Q({}^{\sigma}\!\varphi_{\pi})$ is well-defined up to multiplication by elements of the form $\sigma(t)$ with $t\in E(\pi)^{\times}$ and (within the respective quotient of algebras) independent of the choice of ${}^{\sigma}\!\varphi_{\pi}$. Therefore, the family of numbers $Q({}^{\sigma}\!\varphi_{\pi})$ gives rise to an element
$$Q(\pi)\in E(\pi)\otimes_\Q\C\cong \C^{|J_{E(\pi)}|},$$
called the {\it automorphic $Q$-period} attached to $\pi$.   

Since the Petersson inner product appears in the definition of the GGP-period it is convenient to use it to define the automorphic $Q$-periods. However, it has already been mentioned that the complex conjugation used to define the Petersson inner product does not quite correspond to the operator $F_\infty$ used to define motivic $Q$-periods. Thus, $Q(\pi)$ differs from the $Q$-period attached to the motive $M(\Pi)$
by a factor corresponding to the central character of $\Pi$; this explains the normalization in Definition \ref{local period definition}.

\section{Proof of the factorization}\label{sect:proofThm61}

\subsection{A theorem on critical values of Asai $L$-functions}

\begin{thm}\label{Asai thm}
Let $\Pi$ be a cohomological conjugate self-dual cuspidal automorphic representation of $G_n(\A_F)$, which satisfies Hypothesis \ref{descent}. We assume that $\Pi$ is $5$-regular.  
Then, one has
\begin{equation}\label{Asai L-value 1}
L^{S}(1,\Pi,{\rm As}^{(-1)^{n}})\sim_{E(\Pi)}(2\pi i)^{dn(n+1)/2}\prod\limits_{\imath \in \Sigma}\prod\limits_{0\leq i\leq n}P^{(i)}(\Pi,\imath).
\end{equation}
Interpreted as families, this relation is equivariant under the action of ${\rm Aut}(\C/F^{Gal})$.
\end{thm}

\begin{proof}
By \cite{grob_lin}, Theorem 1.42 \& Theorem 4.17, we know that there is a certain Whittaker-period $p(\Pi)$ attached to $\Pi$ (cf.\ \cite{grob_lin}, Corollary 1.22 \& \S1.5.3), such that  
\begin{equation}\label{eq:q1}
L^{S}(1,\Pi,{\rm As}^{(-1)^{n}})\ \sim_{E(\Pi)} (2\pi i)^{dn} \ p(\Pi).
\end{equation}
 Combining Theorem \ref{thm:artihap} and \cite{jie-thesis}, Corollary $7.5.1$, there exists an archimedean factor $Z(\Pi_{\infty})$ such that 
\begin{equation}\label{eq:q2}
p(\Pi)\sim_{E(\Pi)} Z(\Pi_{\infty})\prod\limits_{\imath\in \Sigma}\prod\limits_{1\leq i\leq n-1}P^{(i)}(\Pi,\imath)\sim_{E(\Pi)} Z(\Pi_{\infty})\prod\limits_{\imath\in \Sigma}\prod\limits_{0\leq i\leq n}P^{(i)}(\Pi,\imath)
\end{equation}
where the last equation follows from equation \eqref{local end 2}.\\\\
Now, let $\Pi^{\#}$ be any cohomological conjugate self-dual cuspidal automorphic representation of $G_{n-1}(\Acm)$, which satisfies Hypothesis \ref{descent} and assume that the coefficient modules in cohomology attached to the pair $(\Pi,\Pi^{\#})$ satisfy the branching-law
\begin{equation}\label{branching law}
\mu_{v,1}\geq -\mu'_{v,n-1}\geq \mu_{v,2}\geq -\mu'_{v,n-2}\cdots \geq -\mu'_{v,1}\geq  \mu_{v,n} 
\end{equation}
for each $v\in S_\infty$. Since $\Pi$ is $5$-regular, we may even choose $\Pi^{\#}$ such that it is also at least $5$-regular. Indeed, such a $\Pi^{\#}$ can be constructed by automorphic induction of an appropriate Hecke character, as in \cite{jie-thesis}, Lemma 9.3.1. The local condition needed for Hypothesis \ref{descent} (see Remark \ref{rmk:des}) can be guaranteed by the argument used in the proof of Theorem \ref{BS}.\\\\  
It follows from \S\ref{sect:critRS}, that there is a critical point $s_0=\tfrac{1}{2}+m$ of $L(s,\Pi\times\Pi^{\#})$ with $m\geq 1$. By definition we have $$Z(\Pi_{\infty})=(2\pi i)^{d(m+\frac{1}{2})n(n-1)-\frac{d(n-1)(n-2)}{2}}\Omega(\Pi_{\infty})^{-1}p(m,\Pi_{\infty},\Pi^{\#}_{\infty})^{-1}$$ (see $(7.8)$ of  \cite{jie-thesis}) where $p(m,\Pi_{\infty},\Pi_{\infty}^{\#})$ is the bottom-degree archimedean Whittaker period attached to $s_0=\tfrac{1}{2}+m$ (see Theorem 1.45 of \cite{grob_lin}) and $\Omega(\Pi^{\#}_{\infty})$ is an archimedean period defined as the ratio of the Whittaker period of an isobaric sum of the product of Whittaker periods of the isobaric summands (see Proposition $3.4.1$ of \cite{jie-thesis} and the last paragraph before section 3.5 of the \textit{loc.cit}).\\\\
One main theorem of \cite{grob_lin} (see Theorem $2.6$ of the \textit{loc.cit}) is that we may define Whittaker periods uniformly such that $\Omega(\Pi_{\infty})\sim 1$. Another main result of the \textit{loc.cit} (see Corollary $4.30$ there) is that 
$$p(m,\Pi_{\infty},\Pi_{\infty}^{\#})\sim_{E(\Pi)E(\Pi^{\#})} (2\pi i)^{mdn(n-1)-\tfrac{1}{2}d(n-1)(n-2)}.
$$
We obtain immediately that $Z(\Pi_\infty)\sim_{E(\Pi)} (2\pi i)^{\tfrac{1}{2}dn(n-1)}$ as claimed.

\end{proof}

\subsection{Statement of the main theorem on factorization}\label{sect:thm}
We shall resume the notation from \S\ref{reviewH14}. In particular, we assume to have fixed a real embedding $\imath_{v_0}$ of $F^+$ and denote by $H=H_{I_0}$ the attached unitary group. Given a highest weight $\lambda$, we obtained $n$ cohomological discrete series representations $\pi_{\lambda,q}$, $0 \leq q \leq n-1$ of $H_\infty$, which were distinguished by the property that their $(\q,K_{H,\infty})$-cohomology is concentrated in degree $q$.\\\\
Now, let $\Pi$ be a cohomological conjugate self-dual cuspidal automorphic representation of $G_n(\A_F)$, which satisfies Hypothesis \ref{descent}. For the same reason as in \S\ref{construction of motive}, we shall descend $\Pi^{{\sf v}}$ instead of $\Pi$. So, for each $q$ as above, we are given a cohomological tempered cuspidal automorphic representation $\pi(q)\in \prod(H,\Pi^{\sf v})$ with archimedean component $\pi_{\lambda,q}$. By Proposition \ref{h,K-coh} it has multiplicity one in the square-integrable automorphic spectrum. Finally, recall the number field $E(\pi(q))\supseteq E_Y(\eta)$ from \S\ref{sect:E(pi)} and let us abbreviate $E_q(\Pi):=E(\Pi) E(\pi(q))$. We are now ready to state our  main result on factorization:

\begin{thm}\label{auto-facto} 
Let $n\geq 2$ and let $\Pi$ be a cohomological conjugate self-dual cuspidal automorphic representation of $G_n(\A_F)$, which satisfies Hypothesis \ref{descent} and let ${\xi}_{\Pi}$ be its central character. We assume that $\Pi_\infty$ is $(n+4)$-regular. Let $\pi(q)\in \prod(H,\Pi^{\sf v})$ be a cohomological tempered cuspidal automorphic representation with archimedean component $\pi_{\lambda,q}$. Moreover we suppose 
Conjecture \ref{lvarch} (rationality of archimedean integrals). Then, for each $0\leq q\leq n-2$, 
$$Q(\pi(q)) \sim_{E_q(\Pi)}  p(\widecheck{\xi}_{\Pi},\Sigma)^{-1} \cfrac{P^{(q+1)}(\Pi,\imath_{v_0}) }{ P^{(q)}(\Pi,\imath_{v_0})}.$$
Interpreted as families, this relation is equivariant under the action of ${\rm Aut}(\C/F^{Gal})$. 
\end{thm}
\noindent Before we give a proof of Theorem \ref{auto-facto}, let us make several remarks and derive an important consequence: \\\\
Firstly, this theorem establishes a version of the factorization of periods which was conjectured in \cite{H97}, see Conjecture 2.8.3 and Corollary 2.8.5 {\it loc.\ cit.}. A proof of this conjecture (up to an unspecified product of archimedean factors) when $\cm=\mathcal K$ is imaginary quadratic was obtained in \cite{H07}, based on an elaborate argument involving the theta correspondence and under a certain regularity hypothesis. The more general argument, which we will give here, is much shorter and more efficient (but evidently depends on the hypotheses of Theorem \ref{auto-facto}).\\\\ 
Secondly, our theorem will imply the desired factorization, cf.\ \eqref{eq:factor}, of the local arithmetic automorphic periods $P^{(i)}(\Pi,\imath)$ as follows:

\begin{defn}\label{local period definition}
Let $\Pi$ and $\pi(q)$ be as in the previous theorem.  We define
$$ P_{i}(\Pi,\imath):= \left\{ \begin{array}{rcl}
         P^{(0)}(\Pi,\imath)& \mbox{if} & i=0; \\
         Q(\pi(i-1)) \ p(\widecheck{\xi}_{\Pi},\Sigma) & \mbox{if}
         & 1\leq i\leq n-1; \\  
P^{(n)}(\Pi,\imath)\prod\limits_{i=0}^{n-1}P_{i}(\Pi,\imath)^{-1} & \mbox{if} & i=n.                \end{array}\right.$$
\end{defn}
Moreover, for any $0\leq i\leq n$, let $E^{(i)}(\Pi)$ be the compositum of the number fields $E_0(\Pi)$ and $E_q(\Pi)$, $q\leq i-1$. Here is our main theorem.

\begin{thm}\label{main factorization}  
Under the regularity hypotheses of Theorem \ref{auto-facto}, and assuming Conjecture \ref{lvarch}, the Tate relation \eqref{eq:fundrel} is true. More precisely, we obtain the following factorization
\begin{equation}\label{Pfac}
P^{(i)}(\Pi,\imath) \sim_{E^{(i)}(\Pi)} P_{0}(\Pi,\imath)P_{1}(\Pi,\imath)\cdots P_{i}(\Pi,\imath)
\end{equation}
and in addition for each $i$ and $\imath$
\begin{equation}\label{eq:Tate}
P_{i}(\Pi,\imath)\sim_{E_{i}(\Pi)}Q_{i}(M(\Pi),\imath),
\end{equation}
for the motive $M(\Pi)$ attached to $\Pi$ as constructed in Theorem \ref{thm:clozel}. Interpreted as families, both relations are equivariant under the action of ${\rm Aut}(\C/F^{Gal})$. 

Assuming only the regularity hypotheses, we have the weaker relations
\begin{equation}\label{eq:Tatew}
P^{(i)}(\Pi,\imath) \sim_{E^{(i)}(\Pi)} I_\infty(\Pi) \cdot P_{0}(\Pi,\imath)P_{1}(\Pi,\imath)\cdots P_{i}(\Pi,\imath);  ~~ 
  P_{i}(\Pi,\imath)\sim_{E_{i}(\Pi)} I_\infty(\Pi) \cdot Q_{i}(M(\Pi),\imath),
\end{equation}
where $I_\infty(\Pi) \in (E_{i}(\Pi)\otimes \CC)^\times$ depends only on the archimedean component $\Pi_\infty$ of $\Pi$.
\end{thm}
\begin{proof}
Given Theorem \ref{auto-facto}, the factorization \eqref{Pfac} follows directly. We now prove \eqref{eq:Tate}: For $i=0$, by equation (\ref{local end}) we have $P^{(0)}(\Pi,\imath)\sim_{E(\Pi)} p(\widecheck{\xi}_{\Pi},\bar\imath) \sim_{E(\Pi)} p(\widecheck{\xi_{\Pi}^{c}},\imath)$. By Definition $3.1$ of \cite{harris-lin}, Equation (2.12) of \cite{linorsay} and Equation (6.13) of \cite{jie-thesis}, we know $$Q_{0}(M(\Pi),\imath) \sim_{E(\Pi)} (2\pi i)^{n(n-1)/2}\delta(M(\Pi),\imath)\sim_{E(\Pi)}\delta(M(\xi_{\Pi}),\imath)\sim_{E(\Pi)}  p(\widecheck{\xi_{\Pi}^{c}},\imath)$$ as expected.\\\\
For each $1\leq i\leq n-1$, by Remark $3.5$ of \cite{H21} (see also Remark \ref{Finfty}), we know 
$$Q_{i}(M(\Pi),\imath) \sim_{E_{i}(\Pi)} Q(\pi(i-1)) q(M(\Pi)),$$ 
where $q(M(\Pi))$ is the period defined in Lemma $4.9$ of \cite{grob-harr}. We see immediately from this lemma that $q(M(\Pi))\sim_{E(\Pi)}\prod\limits_{\imath\in \Sigma}((2\pi i)^{n(n-1)/2}\delta(M(\Pi),\imath))^{-1}$. Similarly as above we have 
$$ (2\pi i)^{n(n-1)/2}\delta(M(\Pi),\imath)\sim_{E(\Pi)}  p(\widecheck{\xi_{\Pi}^{c}},\imath).$$ 
Hence $q(M(\Pi))\sim_{E(\Pi)} \prod\limits_{\imath\in \Sigma} p(\widecheck{\xi_{\Pi}^{c}},\imath)^{-1}\sim_{E(\Pi)}  p(\widecheck{\xi_{\Pi}},\Sigma)$ and $Q_{i}(M(\Pi),\imath) \sim_{E_{i}(\Pi)} P_{i}(\Pi,\imath)$ as expected.\\\\
It remains to show that $\prod\limits_{i=0}^{n}Q_{i}(M(\Pi),\imath)\sim_{E(\Pi)} P^{(n)}(\Pi,\imath)$.  By Lemma $1.2.7$ of \cite{harrisANT} we have $\prod\limits_{i=1}^{n}Q_{i}(M(\Pi),\imath)\sim_{E(\Pi)}((2\pi i)^{n(n-1)/2}\delta(M(\Pi),\imath))^{-2}$. Hence, 
$$\prod\limits_{i=0}^{n}Q_{i}(M(\Pi),\imath)\sim_{E(\Pi)}((2\pi i)^{n(n-1)/2}\delta(M(\Pi),\imath))^{-1}\sim_{E(\Pi)}  p(\widecheck{\xi_{\Pi}^{c}},\imath)^{-1}\sim_{E(\Pi)} p(\widecheck{\xi_{\Pi}},\imath),$$ 
which is equivalent to $P^{(n)}(\Pi,\imath)$ by \eqref{local end}.

\end{proof}


\subsection{Proof of Theorem \ref{auto-facto}}
Our proof will proceed in several steps.

\subsubsection*{Step 1:}
Let us start off with the following 

\begin{obsv}\label{q=0}
Recall that $\pi_{\lambda,q}$ is holomorphic when $q=0$, cf.\ Lemma \ref{knownknown}. It thus follows directly from the definition that we have $Q(\pi(0))\sim_{E(\Pi) E(\pi(0))} P^{(I_{0})}(\Pi)$. Hence, by Theorem \ref{thm:artihap}, and Lemma \ref{Lemma CM},
\begin{eqnarray}
Q(\pi(0))&\sim_{E(\Pi) E(\pi(0))}& \left( \prod\limits_{\imath_v\neq \imath_{v_0}}P^{(0)}(\Pi,\imath_v)\right) P^{(1)}(\Pi,\imath_{v_0}) \nonumber \\
&\sim_{E(\Pi) E(\pi(0)) E_F(\widecheck{\xi}_{\Pi})}& \left( \prod\limits_{\imath_v\neq \imath_{v_0}}p(\widecheck{\xi}_{\Pi},\imath_v)^{-1}\right) P^{(1)}(\Pi,\imath_{v_0}) \nonumber\\
&\sim_{E(\Pi) E(\pi(0)) E_F(\widecheck{\xi}_{\Pi})}&\left( \prod\limits_{\imath_v\in \Sigma}p(\widecheck{\xi}_{\Pi},\imath_v)^{-1}\right) P^{(1)}(\Pi,\imath_{v_0})p(\widecheck{\xi},\bar{\imath}_{v_0})^{-1}\nonumber\\
&\sim_{E(\Pi) E(\pi(0)) E_F(\widecheck{\xi}_{\Pi})}& p(\widecheck{\xi}_{\Pi},\Sigma)^{-1}\cfrac{P^{(1)}(\Pi,\imath_{v_0})}{P^{(0)}(\Pi,\imath_{v_0})}.
\end{eqnarray}
However, by Lemma 1.34 in \cite{grob_lin}, we may reduce this relation to the smallest field containing $F^{Gal}$, on which all the quantities on both sides depend and remain well-defined. But this field is $E_0(\Pi)=E(\Pi) E(\pi(0))$. Therefore, Theorem \ref{auto-facto} is true when $q=0$. This is going to be used as the first step in our inductive argument.
\end{obsv}
Now, let $q$ be arbitrary. Then, in the notation of \S\ref{reviewH14}, the infinity type of $\Pi$ at $v$ is $\{z^{a_{v,i}} \bar{z}^{-a_{v,i}}\}_{1\leq i\leq n}$ where $a_{v,i}=-A_{v,n+1-i}$ (recall that we descend from $\Pi^{\sf v}$ rather than from $\Pi$ now). Next, let $\pi'(q)$ be the cohomological tempered cuspidal automorphic representation of $H'(\A_{F^+})$, constructed in Theorem \ref{BS}. By a direct calculation one gets that the $(\q',K_{H',\infty})$-cohomology of $\pi'(q)$ is non-vanishing only in degree $q':=n-q-2$.\\\\ 
Let $\Pi'=BC(\pi'(q)^{{\sf v}})$ be the base change of the contragredient of $\pi'(q)$. Then, the infinity type of $\Pi'$ at $v\in S_\infty$ is $\{z^{b_{v,j}} \bar{z}^{-b_{v,j}}\}_{1\leq j\leq n-1}$, with $b_{v,j}=A_{v,j}-\tfrac{1}{2}$ if either $v\neq v_{0}$, or $v=v_{0}$ and $j\neq q+1$, whereas $b_{v_{0},q+1}=A_{v_0,q+1}+\tfrac{1}{2}$. Hence, if we calculate the automorphic split indices of the pair $(\Pi,\Pi')$, cf.\ Definition \ref{split automorphic}, then we obtain

$$sp(i,\Pi;\Pi',\imath_v) =     \left\{ \begin{array}{rcl}
         1 & \mbox{if}
         & 1\leq i\leq n-1
 \\ 0  & \mbox{if} & i=0\text{ or }n
                \end{array}\right.;
 sp(j,\Pi';\Pi,\imath_v) =     \begin{array}{rcl}
         1 & \mbox{if}
         & 0\leq j\leq n-1
                \end{array}.
 $$
 at $v\neq v_{0}$, and
 $$sp(i,\Pi;\Pi',\imath_{v_0}) =     \left\{ \begin{array}{rcl}
         1 & \mbox{if}
         & 1\leq i\leq n-1,i\neq q, i\neq q+1
 \\ 0  & \mbox{if} & i=0, q+1 \text{ or }n\\
 2 & \mbox{if} & i=q
                \end{array}\right.;$$
$$sp(j,\Pi';\Pi,\imath_{v_0}) =   \left\{  \begin{array}{rcl}
         1 & \mbox{if}
         & 0\leq j\leq n-1, j\neq n-q-1, j\neq n-q-2\\
2 & \mbox{if}
         & j=n-q-2\\
0 & \mbox{if}
         & j=n-q-1
                \end{array}.\right.
 $$
at $v=v_{0}$. \\\\
We want to insert them into the formula provided by Theorem \ref{automorphic Deligne general intro}: As a first and obvious observation, it is clear by construction that, since $\Pi_\infty$ is $(n+4)$-regular, $\Pi'_\infty$ is $(n+3)$-regular. Combining Theorem \ref{BS}.(1) with Remark \ref{rmk:des}, we also see that $\Pi'$ is a cuspidal automorphic representation, which satisfies Hypothesis \ref{descent}. Therefore, $\Pi'$ satisfies the assumptions of Theorem \ref{automorphic Deligne general intro}, as $n\geq 2$. Hence, inserting the values of the automorphic split indices from above into \eqref{eq:Thm1l}, we obtain

\begin{eqnarray}\label{RS L-value}
\nonumber
L^S(\tfrac{1}{2},\Pi\otimes \Pi')& \sim_{E(\Pi)E(\Pi')} & (2\pi i)^{dn(n-1)/2} \prod\limits_{\imath_v\in \Sigma}\left(\prod\limits_{1\leq i\leq n-1}P^{(i)}(\Pi,\imath_v) \prod\limits_{0\leq j\leq n-1}P^{(j)}(\Pi',\imath_v)\right)\times   \\ &&\cfrac{  P^{(q)}(\Pi,\imath_{v_0})P^{(n-q-2)}(\Pi',\imath_{v_0})    }{  P^{(q+1)}(\Pi,\imath_{v_0}) P^{(n-q-1)}(\Pi',\imath_{v_0})}.
 \end{eqnarray}
 The following observation is crucial for what follows: 
 
\begin{obsv}\label{q1}
$L^S(\tfrac{1}{2},\Pi\otimes \Pi')\neq 0$.
\end{obsv}
\noindent In order to see this, recall that by Theorem \ref{BS}.(2) there are factorizable cuspidal automorphic forms $f \in \pi(q)$, $f' \in \pi'(q)$, whose attached GGP-period does not vanish $\CP(f,f')\neq 0$. Hence, as all the local pairings $I^*_v(f_v,f'_v)$, cf.\ \S\ref{sect:GGPunit}, are convergent by the temperedness of $\pi(q)_v$ and $\pi'(q)_v$, it follows from the Ichino-Ikeda-N.Harris formula, Theorem \ref{conjecture II}, that necessarily
$$L^S(\tfrac{1}{2},BC(\pi(q))\otimes BC(\pi'(q)))\neq 0.$$
But since both $\Pi$ and $\Pi'$ are conjugate self-dual, we have
\begin{equation}\label{eq:l}
L^S(\tfrac{1}{2},BC(\pi(q))\otimes BC(\pi'(q)))=L^S(\tfrac{1}{2},\Pi^{{\sf v}}\otimes \Pi'^{{\sf v}}) =L^S(\tfrac{1}{2},\Pi^{c}\otimes \Pi'^{c}) = L^S(\tfrac{1}{2},\Pi\otimes \Pi').
\end{equation}
Therefore, indeed
$$L^S(\tfrac{1}{2},\Pi\otimes \Pi')\neq 0.$$

\subsubsection*{Step 2:}\label{analysis} 
We resume the notation from Step 1. Recall from the discussion below Theorem \ref{BS} that the factorizable cuspidal automorphic forms $f \in \pi(q)$, $f' \in \pi'(q)$ may be chosen such that, for all $v\in S_\infty$, $f_v$ (resp. $f'_v$) belongs to the $E(\pi(q))$- (resp. $E(\pi'(q))$-) rational subspaces of the minimal $K_{H,v}$-type of $\pi(q)_v$ (resp. $K_{H',v}$-type of $\pi'(q)_v$).\\\\
As in \S\ref{sect:GGPunit}, let $\xi$ be the Hecke character of $U(V_{1})(\A_{F^+})$ given by $\xi=(\xi_{ \pi'(q)} \xi_{\pi(q)})^{-1}$ and write $\pi''(q)=\pi'(q)\otimes \xi$. Let $f_{0}$ be a deRham-rational element of $\xi$. We define $f''=f'\otimes f_{0}$, a deRham-rational element in $\pi''(q)$.  Then, by Lemma \ref{add U1}, the GGP-period 
$$\CP(f,f'')=\cfrac{|I^{can}(f,f'')|^2}{\<f,f\>\ \<f'',f''\>}$$
satisfies
$$\CP(f,f')=\CP(f,f'').$$
Furthermore, by Theorem \ref{BS} and Theorem \ref{GPknownknown}, $I^{can}(f,f'')$ is a non-zero element of $E(\pi(q))E( \pi'(q))$. So, by our choice of $f$ and $f'$, Theorem \ref{conjecture II} and the very definition of the automorphic $Q$-periods attached to $\pi(q)$ and $\pi''(q)$, cf.\ \S\ref{Qperiod}, imply that
\begin{equation}\label{II1}
\begin{aligned}
\frac{1}{Q(\pi(q))Q(\pi''(q))} &\sim_{E(\pi(q)) E(\pi'(q))}   \Delta_{H} \frac{L^S(\tfrac{1}{2},\Pi^{\sf v}\otimes \Pi'^{\sf v})}{L^S(1,\Pi^{\sf v},{\rm As}^{(-1)^n})L^S(1,\Pi'^{\sf v},{\rm As}^{(-1)^{n-1}})} \ \prod_{v \in S_\infty} I^*_v(f_v,f'_v)  \\
&\sim_{E(\pi(q)) E(\pi'(q))}   (2\pi i)^{dn(n+1)/2} \frac{L^S(\tfrac{1}{2},\Pi\otimes \Pi')}{L^S(1,\Pi,{\rm As}^{(-1)^n})L^S(1,\Pi',{\rm As}^{(-1)^{n-1}})} \ \prod_{v \in S_\infty} I^*_v(f_v,f'_v).
\end{aligned}
\end{equation}
Here, we could remove the contragredient in the second line, as both $\Pi$ and $\Pi'$ are conjugate self-dual, whereas the replacement of $\Delta_H$ by a power of $2 \pi i$ is a consequence of (1.37) and (1.38) in \cite{grob_lin}, and the elimination of the local factors $I^*_v(f_v,f'_v)$ at the non-archimedean places follows from Lemma \ref{localIv}. At the archimedean places we make the following observation:

\begin{prop}\label{arch}  
Under the hypotheses of Theorem \ref{BS}, the local factors $I^*_v(f_v,f'_v) \neq 0$ for $v \in S_{\infty}$.
\end{prop}
\begin{proof}  This is an immediate consequence of the non-vanishing of the global period $\CP(f,f')$.
\end{proof}
\noindent Hence, as we are admitting Conjecture \ref{lvarch}, we obtain
\begin{equation}\label{II1b}
\begin{aligned}
\frac{1}{Q(\pi(q))Q(\pi''(q))} &\sim_{E(\pi(q)) E(\pi'(q))}   (2\pi i)^{dn(n+1)/2} \frac{L^S(\tfrac{1}{2},\Pi\otimes \Pi')}{L^S(1,\Pi,{\rm As}^{(-1)^n})L^S(1,\Pi',{\rm As}^{(-1)^{n-1}})}.\end{aligned}
\end{equation}

\subsubsection*{Step 3:}
We recall from Step 1 above that $\Pi'$ is a cohomological conjugate self-dual cuspidal automorphic representation of $G_{n-1}(\A_F)$, which satisfies Hypothesis \ref{descent} and is $(n+3)$-regular. Hence, both $\Pi$ and $\Pi'$ satisfy the conditions of Theorem \ref{thm:artihap} and Theorem \ref{Asai thm}. As a consequence, combining the relations \eqref{local end 2}, \eqref{Asai L-value 1} and \eqref{RS L-value}, one gets 

\begin{equation}\label{II2}
(2\pi i)^{dn(n+1)/2}\frac{L^S(\tfrac{1}{2},\Pi\otimes \Pi')}{L^S(1,\Pi,{\rm As}^{(-1)^n})L^S(1,\Pi',{\rm As}^{(-1)^{n-1}})} \sim_{E(\Pi)E(\Pi')}\cfrac{  P^{(q)}(\Pi,\imath_{v_0})P^{(n-q-2)}(\Pi',\imath_{v_0})    }{  P^{(q+1)}(\Pi,\imath_{v_0}) P^{(n-q-1)}(\Pi',\imath_{v_0})}.
\end{equation}
Recall that $L^S(\tfrac{1}{2},\Pi\otimes \Pi')\neq 0$, cf.\ Observation \ref{q1}. This allows us to combine \eqref{II1b} with \eqref{II2}, and so, using the fact that $Q(\pi''(q)) \sim_{E(\pi'(q))} Q(\pi'(q)) \cdot Q(\xi)$, we arrive at the following conclusion:

\begin{equation}\label{II5}
\frac{1}{Q(\pi(q))Q(\pi'(q))Q(\xi)} \sim_{E_q(\Pi) E_{q'}(\Pi')}   \  \cfrac{  P^{(q)}(\Pi,\imath_{v_0})P^{(n-q-2)}(\Pi',\imath_{v_0})    }{  P^{(q+1)}(\Pi,\imath_{v_0}) P^{(n-q-1)}(\Pi',\imath_{v_0})}
\end{equation}

\subsubsection*{Step 4:}

We need one last ingredient before we can complete the proof of Theorem \ref{auto-facto} by induction on the $F$-rank $n$:

\begin{lem} The following relation $$Q(\xi) \sim_{E_F(\xi)E_F(\widecheck{\xi}_{\Pi})E_F(\widecheck{\xi}_{\Pi'})}    p(\widecheck{\xi}_{\Pi},\Sigma)p(\widecheck{\xi}_{\Pi'}, \Sigma)$$
holds. Interpreted as families of complex numbers it is equivariant under ${\rm Aut}(\C)$.
\end{lem}

\begin{proof}

Recall that $U(V_{1})$ is the one-dimensional unitary group of signature $(1,0)$ at each $\imath\in \Sigma$. Let $T_1:=R_{F/\Q}(U(V_1))$. By definition of the CM-periods we have $Q(\xi)\sim_{E_F(\xi)} p(\xi, (T_1,h_{1}))$ where $h_{1}: R_{\C/\R}(\mathbb{G}_{m,\C})\rightarrow T_{1,\R}$ is the map, which sends $z$ to $z/\overline{z}$ at each $\imath\in\Sigma$. \\\\
We define a map $h_{\tilde{\Sigma}}: R_{\C/\R}(\mathbb{G}_{m,\C})\rightarrow T_{F,\R}$, where $T_{F}=R_{F/\Q}(\mathbb{G}_{m})$, by sending $z$ to $z/\overline{z}$ at each $\imath\in \Sigma$. The pair $(T_{F},h_{\tilde{\Sigma}})$ is then a Shimura datum. We extend $\xi$ to a character of $\Acm^{\times}$, still denoted by $\xi$. The natural inclusion $T_1\hookrightarrow T_{F}$ induces a map from the Shimura datum $(T_1,h_1)$ to $(T_{F},h_{\tilde{\Sigma}})$. By Proposition \ref{relation CM-period}, we have $p(\xi, (T_1,h_1))\sim_{E_F(\xi)} p(\xi, (T_{F},h_{\tilde{\Sigma}}))$.\\\\
Let $(T_{F},h_{\Sigma})$ and $(T_{F},h_{\overline{\Sigma}})$ be as in \S\ref{CM-periods}. Multiplication defines a map from $(T_{F}, h_{\tilde{\Sigma}})\times (T_{F},h_{\overline{\Sigma}})$ to $(T_{F}, h_{\Sigma})$. It follows from Proposition \ref{relation CM-period} (see also Proposition $1.4$ and Corollary $1.5$ of \cite{Harris93}), that we have 
$$
p(\xi, (T_{F},h_{\tilde{\Sigma}}))\sim_{E_F(\xi)} p(\xi,(T_{F},h_{\Sigma}))p(\xi,(T_{F},h_{\overline{\Sigma}}))^{-1}=p(\xi,\Sigma)p(\xi, \overline{\Sigma})^{-1}.$$
By Lemma \ref{Lemma CM},
$p(\xi,\Sigma)p(\xi, \overline{\Sigma})^{-1}\sim_{E_F(\xi)} p(\xi,\Sigma)p(\xi^{c,-1}, \Sigma)\sim_{E_F(\xi)} p(\xi/\xi^{c},\Sigma)$. Note that $\xi/\xi^{c}$ is the base change of the original $\xi$. Recall that $\Pi^{c}\cong \Pi^{{\sf v}}$ is the base change of $\pi(q)$. Hence $\widecheck{\xi}_{\Pi}$ is the base change of $\xi_{\pi(q)}^{-1}$. Similarly, $\widecheck{\xi}_{\Pi'}$ is the base change of $\xi_{\pi'(q)}^{-1}$. Therefore, $\xi/\xi^{c}=\widecheck{\xi}_{\Pi}\widecheck{\xi}_{\Pi}$. Consequently, recollecting all relations from above and invoking Lemma \ref{Lemma CM} once more, we get
\begin{equation}
Q(\xi)\sim_{E_F(\xi)} p(\widecheck{\xi}_{\Pi}\widecheck{\xi}_{\Pi'},\Sigma)\sim_{E_F(\xi)E_F(\widecheck{\xi}_{\Pi})E_F(\widecheck{\xi}_{\Pi'})} p(\widecheck{\xi}_{\Pi},\Sigma)p(\widecheck{\xi}_{\Pi'},\Sigma).
\end{equation}
\end{proof}
The previous Lemma and equation (\ref{II5}) now implies

\begin{equation}\label{II6}
Q(\pi(q))Q(\pi'(q)) \sim_{E_q(\Pi) E_{q'}(\Pi')}   \   \left( p(\widecheck{\xi}_{\Pi},\Sigma)^{-1}\cfrac{ P^{(q+1)}(\Pi,\imath_{v_0})}{P^{(q)}(\Pi,\imath_{v_0})}\right) \times  \left( p(\widecheck{\xi}_{\Pi'},\Sigma)^{-1} \cfrac{P^{(n-q-1)}(\Pi',\imath_{v_0})    }{   P^{(n-q-2)}(\Pi',\imath_{v_0})}\right)
\end{equation}
Here, we could remove the number field  $E_F(\xi)E_F(\widecheck{\xi}_{\Pi})E_F(\widecheck{\xi}_{\Pi'})$ from the relation using \cite{grob_lin}, Lemma 1.34.\\\\
We may finish the proof of Theorem \ref{auto-facto} by induction on $n$. When $n=2$, the integer $q$ is necessarily $0$. The theorem is then clear by Observation \ref{q=0}. We assume that the theorem is true for $n-1\geq 2$. Again, if $q=0$, then the theorem follows from Observation \ref{q=0}. So, let $1\leq q\leq n-2$. Recall that our representation $\pi'(q)$ from above is an element in $ \prod(H',\Pi'^{{\sf v}})$ whose $(\q',K_{H',\infty})$-cohomology in concentrated in degree $n-q-2\leq n-3$. Moreover, we have verified above that $\Pi'$ is $(n-2)$-regular and satisifes Hypothesis \ref{descent}, whence $\Pi'$ and $\pi'(q)$ satisfy the conditions of Theorem \ref{auto-facto}. Hence $Q(\pi'(q)) \sim_{E_{q'}(\Pi')}   \    p(\widecheck{\xi}_{\Pi'},\Sigma)^{-1} \cfrac{P^{(n-q-1)}(\Pi',\imath_{v_0})    }{   P^{(n-q-2)}(\Pi',\imath_{v_0})}$. The theorem then follows from equation (\ref{II6}) and \cite{grob_lin}, Lemma 1.34.

\subsection{Refinement of Theorem \ref{automorphic Deligne general intro}}

We conclude by restating Theorem \ref{cplusintro} explicitly as a refinement of Theorem \ref{automorphic Deligne general intro}.

\begin{thm}\label{mainDeligneconj}
Let $n,n'\geq 1$ be integers and let $\Pi$ (resp. $\Pi'$) be a cohomological conjugate self-dual cuspidal automorphic representation of $G_n(\Acm)$ (resp. $G_{n'}(\Acm)$), which descends to a tempered cuspidal automorphic representation of a unitary group $U_I(\A_{F^+})$ for each possible $[F^+:\Q]$-tuple of signatures $I$ at the archimedean places, i.e., it satisfies Hypotheis \ref{descent}. We assume that both $\Pi_\infty$ and $\Pi'_{\infty}$ are $5$-regular.  If $n \equiv n' \mod 2$, we assume in addition that the isobaric sum $(\Pi\eta^n)\boxplus (\Pi'^{c}\eta^{n'})$ is $2$-regular; if $n$ and $n'$ have opposite parities then we assume $(\Pi\eta^n)\boxplus (\Pi'^{c}\eta^{n'})$ is $5$-regular.  Then the following version of Deligne's conjecture, cf.\ Conjecture \ref{main conjecture}, is true: If $s_0$ is critical, in Deligne's sense, for $L(s,\Pi\times \Pi')$, then the value at $s_0$ of the partial $L$-function $L^S(s,\Pi\times \Pi')$ (for some appropriate finite set $S$), satisfies
\begin{equation}\label{eq:Thmfinal}
L^S(s_{0},\Pi\otimes \Pi') \sim_{E(\Pi)E(\Pi')} I_{\infty}(\Pi,\Pi')c^+(s_0,R_{F/\Q}(M(\Pi)\otimes M(\Pi'))).
\end{equation}
Here $I_\infty(\Pi,\Pi') \in (E(\Pi)E(\Pi')\otimes \CC)^\times$ depends only on the archimedean factors of $\Pi$ and $\Pi'$, and $c^+(s_0,R_{F/\Q}(M(\Pi)\otimes M(\Pi')))$ is the version of Deligne's period, as recalled in \S \ref{general Deligne} obtained from the motives $M(\Pi)$ and $M(\Pi')$ constructed in Theorem \ref{thm:clozel}.

Using Proposition \ref{Deligne period motivic} to the motives $M(\Pi)$ and $M(\Pi')$, and letting $\mathbf{d}(n,n') = \frac{nn'd(n+n'-2)}{2}$, we thus have
\begin{eqnarray}
&L^S(s_{0},\Pi\otimes \Pi') \sim_{E(\Pi)E(\Pi')} &\\ \nonumber
 &I_{\infty}(\Pi,\Pi') (2\pi i)^{-\mathbf{d}(n,n')} \prod\limits_{\imath\in \Sigma} [\prod\limits_{j=0}^{n}Q^{(j)}(M(\Pi),\imath)^{sp(j,M(\Pi);M(\Pi'),\imath)}\prod\limits_{k=0}^{n'}Q^{(k)}(M(\Pi'),\imath)^{sp(k,M(\Pi');M(\Pi),\imath)}],
 \end{eqnarray}
 with the invariants $Q^{(\bullet)}$ defined by  the relations \eqref{eq:Qperiods} and Theorem \ref{main factorization}.

\end{thm}

\vskip 10pt

\footnotesize
{\sc Harald Grobner: Fakult\"at f\"ur Mathematik, University of Vienna, Oskar--Morgenstern--Platz 1,  A-1090 Vienna, Austria.}
\\ {\it E-mail address:} {\tt harald.grobner@univie.ac.at}

\vskip 10pt

{\sc Michael Harris: Department of Mathematics, Columbia University, New York, NY  10027, USA. }
\\ {\it E-mail address:} {\tt harris@math.columbia.edu}

\vskip 10pt

{\sc Jie Lin: Universit\"at Duisburg-Essen, Fakult\"at f\"ur Mathematik, Mathematikcarr\'ee, Thea-Leymann-Strasse 9,
45127 Essen, Germany}
\\ {\it E-mail address:} {\tt jie.lin@uni-due.de}

\end{document}